\documentclass[dvipsnames,reqno]{siamart250211}
\newsiamremark{remark}{Remark}
\newsiamremark{hypothesis}{Hypothesis}
\crefname{hypothesis}{Hypothesis}{Hypotheses}
\newsiamthm{claim}{Claim}
\newcommand{\mynewline}{\newline}

\usepackage{booktabs}
\usepackage{tabularx}
\usepackage{mathtools}
\usepackage{amsmath}
\usepackage{amsfonts}
\usepackage{amssymb}
\usepackage{subfigure}
\usepackage{verbatim}
\usepackage{cancel}
\newcommand{\stkout}[1]{\cancel{#1}}

\usepackage{tikz}
\usepackage{ifthen}
\usetikzlibrary{matrix}
\usetikzlibrary{positioning}
\usetikzlibrary{shapes.geometric}
\usepackage{pgfplots}
\pgfplotsset{compat=1.18}
\usepgfplotslibrary{groupplots}

\usepackage[algosection,vlined,algo2e,linesnumbered,algoruled]{algorithm2e}
\DontPrintSemicolon

\usepackage{cleveref}

\crefname{algocf}{algorithm}{algorithms}
\Crefname{algocf}{Algorithm}{Algorithms}

\newcommand{\emin}{\ensuremath{e_\mathrm{min}}}
\newcommand{\emax}{\ensuremath{e_\mathrm{max}}}

\DeclareMathOperator{\diag}{diag}

\DeclarePairedDelimiter{\abs}{\lvert}{\rvert}
\DeclarePairedDelimiter{\ceil}{\lceil}{\rceil}
\DeclarePairedDelimiter{\floor}{\lfloor}{\rfloor}

\newcommand{\numset}[1]{\ensuremath{\mathbb{#1}}}
\newcommand{\N}[1][]{\numset{N}^{#1}}
\newcommand{\R}{\numset{R}}

\newcommand{\FIN}{\ensuremath{\mathbb}{F}\langle \emin, \emax, p\rangle}

\DeclareMathOperator{\fl}{\operatorname{fl}}

\newcommand{\sm}{\ensuremath{s_{m}}}
\newcommand{\sM}{\ensuremath{s_{M}}}
\newcommand{\sA}{\ensuremath{s_{A}^{\vphantom{*}}}}
\newcommand{\sB}{\ensuremath{s_{B}^{\vphantom{*}}}}
\newcommand{\sAs}{\ensuremath{s_{A}^{*}}}
\newcommand{\sBs}{\ensuremath{s_{B}^{*}}}
\newcommand{\uA}{\ensuremath{u_{A}}}
\newcommand{\uB}{\ensuremath{u_{B}}}
\newcommand{\kA}{\ensuremath{\kappa_{A}}}
\newcommand{\kB}{\ensuremath{\kappa_{B}}}
\newcommand{\DA}{\varDelta A}
\newcommand{\DB}{\varDelta B}
\newcommand{\da}{\delta a}
\newcommand{\db}{\delta b}

\newcommand{\bydim}[2]{{\ensuremath{#1 \times #2}}}
\newcommand{\kbyn}{\bydim{k}{n}}
\newcommand{\mbyk}{\bydim{m}{k}}

\newcommand{\mbyn}{\bydim{m}{n}}

\newcommand{\F}{\numset{F}}
\newcommand{\I}{\numset{I}}
\newcommand{\fpf}[1][p]{\ensuremath{\F_{#1}}}
\newcommand{\Fkn}{\fpf^{\kbyn}}
\newcommand{\Fmk}{\fpf^{\mbyk}}
\newcommand{\Fmn}{\fpf^{\mbyn}}

\newcommand{\intf}[1][t]{\ensuremath{\I_{#1}}}
\newcommand{\Imk}{\intf^{\mbyk}}
\newcommand{\Ikn}{\intf^{\kbyn}}

\newcommand{\wh}{\widehat}

\newcommand{\thetitle}{Analysis of Floating-Point Matrix Multiplication Computed via Integer Arithmetic}

\newcommand{\thefunding}{%
The work of the last author was supported by Engineering and Physical Sciences Research Council grant EP/W018101/1.}

\newcommand{\theauthori}{Ahmad Abdelfattah}
\newcommand{\theauthorii}{Jack Dongarra}
\newcommand{\theauthoriii}{Massimiliano Fasi}
\newcommand{\theauthoriv}{\mynewline{}Mantas Mikaitis}
\newcommand{\theauthorv}{Françoise Tisseur}

\newcommand{\theaffiliationi}{%
    Innovative Computing Laboratory,
    University of Tennessee,
    Knoxville, TN, USA
    (\email{ahmad@icl.utk.edu},
    \email{dongarra@icl.utk.edu})}
\newcommand{\theaffiliationii}{%
    Department of Mathematics,
    University of Manchester,
    Oxford Road, Manchester M13 9PL, UK
    (%
    \email{francoise.tisseur@manchester.ac.uk})}
\newcommand{\theaffiliationiii}{%
    School of Computer Science,
    University of Leeds,
    Woodhouse Lane, Leeds LS2 9JT, UK
    (\email{m.fasi@leeds.ac.uk},
    \email{m.mikaitis@leeds.ac.uk})}

\newcommand{\theshortauthor}{A.~Abdelfattah, J.~Dongarra, M.~Fasi, M.~Mikaitis, and F.~Tisseur}

\newcommand{\theabstract}{%
    Ootomo, Ozaki, and Yokota [\textit{Int.~J.~High Perform.~Comput.~Appl.}, 38 (2024), p.~297--313] have proposed a strategy to recast a floating-point matrix multiplication in terms of integer matrix products.
    The factors $A$ and $B$ are split into integer \emph{slices}, the product of these slices is computed exactly, and $AB$ is approximated by accumulating these integer products in floating-point arithmetic.
    This technique is particularly well suited to mixed-precision matrix multiply--accumulate units with integer support, such as the NVIDIA tensor cores or the AMD matrix cores.
    The number of slices allows for performance-accuracy tradeoffs: more slices yield better accuracy but require more multiplications, which in turn reduce performance.
    We propose an inexpensive way to estimate the minimum number of multiplications needed to achieve a prescribed level of accuracy.
    Our error analysis shows that the algorithm may become inaccurate (or inefficient) if rows of $A$ or columns of $B$ are \emph{badly scaled}.
    We perform a range of numerical experiments, both in simulation and on the latest NVIDIA GPUs, that confirm the analysis and illustrate strengths and weaknesses of the algorithm.}

\newcommand{\theshorttitle}{Matrix Multiplication with Integer Arithmetic}

\newcommand{\thekeywords}{
    matrix multiplication,
    floating-point arithmetic,
    integer arithmetic,
    tensor cores,
    mixed-precision,
    fixed-point arithmetic,
    error analysis.
  }

\newcommand{\themsccodes}{%
65F99, %
65G50, %
65Y10 %
}

\makeatletter
\newcommand*{\addFileDependency}[1]{%
  \typeout{(#1)}%
  \@addtofilelist{#1}%
  \IfFileExists{#1}{}{\typeout{No file #1.}}%
}
\makeatother

\pgfplotsset{
  height=1.8in,
  width=2.3in,
}

\pgfplotsset{
    reference style/.style={
        thick,
        Apricot,
        solid,
        mark=*,
        mark options={draw=none}
    },
    algorithm base style/.style={
        thick,
        mark=*,
        mark size=2pt,
        mark options={solid}
    },
    ozimmu style/.style={
        algorithm base style,
    },
    cuimma style/.style={
        algorithm base style,
        BurntOrange,
        mark=triangle*,
    },
    ozimmu s 2 style/.style={ozimmu style, mark=square, solid, ProcessBlue},
    ozimmu s 3 style/.style={ozimmu style, mark=triangle, densely dashed, mark options={rotate=0, solid}, Blue},
    ozimmu s 4 style/.style={ozimmu style, mark=o, densely dashed, RoyalBlue},
    ozimmu s 6 style/.style={ozimmu style, mark=triangle, dashed, mark options={rotate=90, solid}, TealBlue},
    ozimmu s 7 style/.style={ozimmu style, mark=triangle, densely dotted, mark options={rotate=180, solid}, NavyBlue},
    ozimmu s 8 style/.style={ozimmu style, mark=triangle, densely dashdotted, mark options={rotate=270, solid}, Cerulean},
    ozimmu s 10 style/.style={ozimmu style, mark=pentagon, densely dotted, Turquoise},
    cuimma s 3 style/.style={cuimma style, mark=diamond, densely dashed, BrickRed},
    cuimma s 6 style/.style={cuimma style, mark=square, dashed, RawSienna},
    cuimma s 7 style/.style={cuimma style, mark=pentagon, densely dotted, RedOrange},
    cuimma s 8 style/.style={cuimma style, mark=square, densely dashdotted, Sepia},
    cuimma s 12 style/.style={cuimma style, mark=10-pointed star, loosely dashed, VioletRed},
    cuimma s 16 style/.style={cuimma style, mark=x, dotted, Plum},
    cuimma s 18 style/.style={cuimma style, mark=+, loosely dotted, Magenta},
}

\begin{document}

\headers{\theshorttitle}{\theshortauthor}

\title{\thetitle\footnotemark[1]}
\date{}

\author{%
  \theauthori{}\footnotemark[2] \and
  \theauthorii{}\footnotemark[2] \footnotemark[3] \and
  \theauthoriii{}\footnotemark[4] \and
  \theauthoriv{}\footnotemark[4] \and
  \theauthorv{}\footnotemark[3]
}

\maketitle

\renewcommand{\thefootnote}{\fnsymbol{footnote}}
\footnotetext[1]{Version of 16 April 2026. \funding{\thefunding{}}}
\footnotetext[2]{\theaffiliationi}
\footnotetext[3]{\theaffiliationii}
\footnotetext[4]{\theaffiliationiii}

\renewcommand{\thefootnote}{\arabic{footnote}}
\setcounter{footnote}{0}

\begin{abstract}
\theabstract
\end{abstract}

\section{Introduction}

The top four computers on the November 2025 Top500 list,\footnote{\url{https://www.top500.org/lists/top500/list/2025/11/}} El Capitan, Frontier, Aurora, and JUPITER are exascale supercomputers, capable of performing over $10^{18}$ floating-point operations per second (flop/s) in binary64 arithmetic, whose accuracy is essential for most scientific applications.

To achieve their impressive performance, modern supercomputers leverage hardware accelerators designed for machine-learning workloads, which typically do not require high precision and can provide meaningful results if fewer-than-32-bit floating-point arithmetics are used.
Formats such as TensorFloat-32, bfloat16, and binary16 are widely available in hardware, and more recently vendors have started developing 8-bit formats for training and inference of deep neural networks: Graphcore has proposed three such formats~\cite{njjm22}, two of which are available in the Tile Vertex ISA~\cite{grap22}; NVIDIA, Arm, and Intel have proposed two~\cite{msbc22}, subsequently crystallized in the Open Compute Project 8-bit floating-point specification (OFP8)~\cite{modc23}; Tesla has proposed the Configurable Float8 format in its Dojo Technology white paper~\cite{tesl21}; and Huawei has proposed the Ascend HiFloat8 format~\cite{lzwl24}. More examples can be found in the interim report of the IEEE P3109 working group~\cite{p3109_24}, which is currently developing a standard for arithmetic formats for machine learning.
To further complicate the landscape, integer arithmetic is often preferred for inference~\cite{bknr23}, and hardware accelerators are starting to be optimized for compact integer formats.
The main features of these formats are summarised in \cref{sec:integ-float-point}.

These reduced-precision formats can have a throughput over two orders of magnitude higher than binary64, but they lack the precision needed for traditional scientific simulations, which require higher accuracy to deliver meaningful results.
Currently, supercomputer-grade accelerators support binary64 arithmetic, but the field is shifting, and soon many will be optimized exclusively for lower precision, with support for even binary32 arithmetic expected to decline over the coming years.
In fact, the performance gap between high- and low-precision arithmetic is already so wide one \emph{must} rely on low-precision formats to fully utilize these accelerators.
Therefore, to integrate GPUs effectively into scientific computing, computations must be re-imagined to use mixed precision, recasting high-precision operations in terms of low-precision ones, possibly resorting to integer arithmetic whenever feasible.

In the context of numerical linear algebra, mixed-precision algorithms have been developed to compute matrix products~\cite{fhlm23},~\cite{fwcz21},~\cite{hth19},~\cite{lxlw21},~\cite{mwfz22},~\cite{mclp18},~\cite{mooi20},~\cite{mooi21},~\cite{pili21},~\cite{vjlv23}, and to solve efficiently linear systems~\cite{abhl24},~\cite{cahi17},~\cite{cahi18},~\cite{hipr21} and least squares problems~\cite{chp20}.
We refer the reader to the survey by Abdelfattah et al.~\cite{aabc21} for a broad overview of existing mixed-precision numerical algorithms for linear algebra, and to the work of Higham and Mary~\cite{hima22} for a discussion of the underlying error analysis.

Matrix multiplication is arguably one of the most fundamental linear algebra kernels: it underpins the majority of numerical algorithms used in matrix computations---small improvements in its performance can have a significant impact on the overall execution time of many matrix algorithms, including the standard operations provided by the LAPACK interface.
Here we focus on the \emph{Ozaki scheme} for matrix multiplication.
This technique can be traced back to the seminal work of Ozaki, Ogita, Oishi, and Rump~\cite{ooor12},~\cite{ooor13} and is based on a technique for accurate floating-point summation due to Rump, Ogita, and Oishi~\cite{roo08}.
Given two floating-point matrices $A$ of size $m \times k$ and $B$ of size $k \times n$, the algorithm computes the product $AB$ in three steps.
First, the rows of $A$ are converted to a block fixed-point representation, where all values in a row share the same scaling factor, and the columns of $B$ are similarly converted to a fixed-point representation.
These fixed-point representations, which need not be obtained explicitly, are then split into slices, where the number of significant bits in each slice is chosen so that the product of two slices can be computed exactly on the mixed-precision hardware available.
Finally, the exact partial products are accumulated in floating-point arithmetic to yield the final result.

Early work on the Ozaki scheme~\cite{ooor12},~\cite{ooor13} considers binary64 arithmetic not only for the input and output matrices, but also for the intermediate slices.
In 2017, NVIDIA introduced the first GPUs featuring \emph{tensor cores}, mixed-precision units that can multiply matrices of binary16 values and accumulate results using binary32 arithmetic.
Whilst some GPUs are equipped with tensor cores that support binary64 arithmetic, only low precision yields high performance.
But this comes at a cost, because low precision variants produce less accurate results and do not use standard IEEE 754-conforming arithmetic~\cite{fhmp21}.
Mukunoki et al.~\cite{mooi20} designed a version of the Ozaki scheme that used these early tensor cores to multiply matrices of binary32 and binary64 values using binary16 slices.
The same ideas were later applied to the computation of products of matrices of binary128 values using binary64 slices~\cite{mooi21}.

In 2018, the NVIDIA Turing GPU Architecture~\cite{nvid18} introduced the second generation of tensor cores, which supports the multiplication of matrices of 8-bit signed integers, with results stored in a 32-bit signed integer format.
The subsequent performance improvement in the Ampere~\cite{nvid20} and Hopper~\cite{nvid22} microarchitectures led Ootomo, Ozaki, and Yokota~\cite{ooy24} to propose a variant of the Ozaki scheme that uses 8-bit signed integer slices to compute the product of two matrices of binary64 values.
These algorithms are summarized in \cref{sec:ozaki-scheme}.
The accuracy and performance of this technique were later optimized by Uchino, Ozaki, and Imamura~\cite{uoi25}, while Lin et al.~\cite{lszl24} considered the use of integer arithmetic to compute the product of matrices of binary32 values.

Our contribution is two-fold.
First, we propose a new error analysis of the integer Ozaki scheme, which can be found in \cref{sec:error-analysis}.
This new analysis shows that the scheme may fail if $A$ has badly scaled rows or $B$ has badly scaled columns.
By considering a different number of slices for $A$ and $B$, the analysis offers an inexpensive strategy to minimize the number of matrix multiplications while providing a bound on the accuracy of the final result.

Second, we run a range of numerical experiments, whose results are reported in \cref{sec:numerical-experiments}, to assess the accuracy and performance of the Ozaki scheme.
Our analysis suggests that this approach may require an extremely large number of slices if the matrices are badly scaled.
In \cref{sec:minimal-example}, we illustrate this with a minimal hand-crafted example and with large randomly-generated, badly-scaled matrices.
Our analysis also suggests that, in some cases, using a different number of slices for the two input matrices will cause no loss in accuracy.
This is illustrated in \cref{sec:lu-experiment}, where we use a block LU factorization to solve linear systems with matrices from literature.
To understand the performance of the scheme, in \cref{sec:benchmarking-gpus} we compare different high-performance implementations of the Ozaki scheme on modern GPU architectures.
To our knowledge, we are the first to report performance results for this algorithm on an NVIDIA Blackwell GPU.
Our results confirm that the reduced-product variants of the Ozaki scheme work well in practice under favorable scaling conditions and when the slices are chosen judiciously.

\section{Integer and floating-point arithmetic}
\label{sec:integ-float-point}

Any nonzero real number $x \in \R$ can be expressed in normalized scientific notation as
\begin{equation*}
  x=(-1)^s\cdot 2^e \cdot m,
\end{equation*}
where $s \in \{0,1\}$ is the sign, $e \in \mathbb{Z}$ is the exponent, and $m \in [1,2)$ is the significand.
In floating-point arithmetic, we discretize the reals by limiting the maximum precision of $m$ and the range of $e$.
Here, we consider the floating-point number system $\FIN$, which is the finite subset of $\R$ obtained by restricting the maximum number of significant bits in $m$ to $p>0$ and by requiring that $\emin \leq e \leq \emax$, with $\emin < \emax$.
The significand $m \in [1, 2)$ is a real number with at most $p$ binary digits (bits).
The requirement that $m$ be between 1 and 2 is usually relaxed for $e = \emin$, in which case $m$ is allowed to be any $p$-bit positive real no greater than 2.
We use the shorthand notation $\fpf \equiv \FIN$ whenever the values of $\emin$ and $\emax$ are clear from the context.

\newcommand{\expbits}{\ensuremath{b_{e}}}
We assume that the floating-point numbers are encoded as binary strings using the encoding in~\cite[sect.~3.4]{ieee19}.
The sign is stored in the leftmost bit of the representation, so that the number is negative if the sign bit is set and positive otherwise.
The $\expbits$ bits immediately to the right of the sign bit are used to store the exponent using a biased representation.
The IEEE 754 format construction rules require that $\emin=1-\emax$, in which case we can set $\emax = 2^{\expbits-1}-1$ and use a representation biased by $\emax$.
Therefore, the smallest and largest allowed exponents are represented as $00\cdots01_2$ and $11\cdots10_2$, respectively.
The all-zero string is reserved for subnormal numbers, whose exponent is $\emin$, and the all-one string is reserved for special values that are needed to ensure that the semantics of all floating-point operations are well specified.
The remaining bits are used to store the fraction, which contains the trailing $p-1$ bits of the significand of $x$, as the left-most bit can be inferred from the exponent field: it will be a zero if the exponent field is the all-zero string, and a one if the exponent field is neither the all-zero nor the all-one string.

For binary64, which is the IEEE format of interest in this work, $p = 53$ and $\expbits = 11$, which implies $\emax = 1023$.

We will consider the round-to-nearest function $\fl : \R \to \fpf$, which maps a real $\R$ to the closest element of $\fpf$.
Regardless of the rule used to break ties, it can be shown that this rounding function satisfies the property~\cite[Thm.~2.2]{high02}
\begin{equation*}
  \fl(x) = x (1 + \delta), \qquad \abs{\delta} < u,
\end{equation*}
where $u = 2^{-p}$ is the unit roundoff of $\fpf$.
In general, the result of a computation involving numbers in $\fpf$ is not an element of $\fpf$.
We assume that the relative error in the result follows the standard model of floating-point arithmetic~\cite[eq.~(2.4)]{high02}, which states that for any $x, y \in \fpf$, the elementary arithmetic operations satisfy
\begin{equation*}
    \fl(x \circ y) = (x \circ y)(1 + \delta),\qquad \abs{\delta} < u,\qquad \circ \in \{+, -, \times, \div\},
\end{equation*}
and a similar result is usually assumed for square root.

We will denote by $\intf$ a signed integer format that uses $t+1$ bits.
We will assume that the numbers are stored using the two's complement representation, so that $\intf$ can represent integers in $[-2^{t},2^{t}-1]$.
A key point of integer arithmetic is that addition and multiplication are exact unless the result overflows.
If two's complement is used, representing the product of two elements in $\intf$ requires at most $2t+1$ bits~\cite[p.~31]{erla04}, and representing their sum will require at most $t+2$ bits~\cite[p.~17]{erla04}.
Adding more than two integers requires a larger number of extra bits, but it is well known~\cite[p.~138]{erla04} that adding up $k$ $t$-bit integers will require at most
\begin{equation}
  \label{eq:max-num-bits-sum-n}
  t + \ceil{\log_{2}k}
\end{equation}
bits.
Fixed-point arithmetic is usually implemented using integers, as the position of the binary point is the only additional information needed to convert a binary integer representation to a fixed-point one.

\paragraph{Matrix multiply--accumulate units}
The simulation of scalar floating-point operations using integer arithmetic is well understood~\cite[sect.~4.2.1]{knut97},~\cite[Chap.~9]{mbdj18}, and efficient software implementations, such as SoftFloat~\cite{haus18} and FLIP~\cite{jere18}, are available.
The algorithms considered here address a different problem: they operate at the matrix level and simulate floating-point matrix multiplication by exploiting integer matrix multiply--accumulate (MMA) operations.
These units can compute $AB + C$ where $A$, $B$, and $C$ are matrices of signed integers represented in two's complement.
The algorithms addressed below also utilize binary64 scalar addition operations.
The NVIDIA Hopper~\cite{nvid22} and Blackwell~\cite{nvid24a} microarchitectures contain tensor cores that support INT8 input datatype.
The NVIDIA PTX ISA 8.7~\cite[sect.~9.7.16]{nvid25a} lists 4- and 8-bit unsigned and signed integer formats for $A$ and $B$ and a 32-bit signed integer format for the accumulators $C$ and $D$.

\section{The Ozaki scheme with integer block MMA}
\label{sec:ozaki-scheme}

The Ozaki scheme~\cite{ooor12,ooor13} is an algorithm for matrix multiplication that exploits the error-free transformation for accurate floating-point summation proposed by Rump, Ogita, and Oishi~\cite{roo08}.
Traditionally, the Ozaki scheme used floating-point arithmetic throughout.
For example, Mukunoki et al.~use it to implement accurate binary32 and binary64 matrix--matrix multiplication on NVIDIA GPUs equipped with first-generation NVIDIA tensor cores~\cite{mooi20}, and to achieve binary128 accuracy using only binary64 matrix multiplication~\cite{mooi21}.
Recently, however, Ootomo, Ozaki, and Yokota~\cite{ooy24} have proposed a variant that relies on integer matrix multiplication and is expected to be very efficient on the upcoming generation of NVIDIA GPUs.
In this section, we review this algorithm and an improved version proposed by Uchino, Ozaki, and Imamura~\cite{uoi25}.

\paragraph{Original algorithm}

Let $A \in \Fmk$ and $B \in \Fkn$ be matrices with no infinities, NaNs, or negative zeros.
In this section and following sections, we further assume that computation does not produce infinities and NaNs.
Rows of $A$ and columns of $B$ with only zeros do not affect the result, thus we assume that each row of $A$ and column of $B$ contains at least one nonzero element.
We describe the algorithm in its full generality, but to aid the reader we provide a small worked example.
In the example, we consider an inner product ($m = n = 1$) with $k = 3$, and for the formats we set $p = 8$ and use the integer format $\intf[t]$ with $t = 3$.
The two vectors we consider are
\begin{equation}
    \label{eq:ex-mats}
  A = \begin{bmatrix} 1.5625    &8    &-3.6875 \end{bmatrix},\qquad
  B = \begin{bmatrix} 1.3828125\\-7.625\\3.625 \end{bmatrix}.
\end{equation}

At a high level, the integer Ozaki scheme of Ootomo, Ozaki, and Yokota~\cite{ooy24} approximates the product $AB$ in three steps.
First, the entries of $A$ and $B$ are implicitly converted to a block fixed-point representation, where a block is a row of $A$ or a column of $B$.
In this format, a block contains a number of \emph{elements}, and all elements in a block share the same \emph{scale}, which is a power of 2.
All entries in the $i$th row of $A$ share a single scale factor, $\alpha_{i}$, defined by
\begin{equation}
  \label{eq:alpha-def}
  \alpha_{i} = 2^{\floor{\log_{2} M_{i}} + 1},\qquad
  M_{i}= \max_{1\le j \le k} \abs{a_{ij}},\qquad
  1 \le i \le m,
\end{equation}
which is the smallest power of two that is larger than the maximum value (in magnitude) within the block.
This guarantees that $0 \le \abs{a_{ij} / \alpha_{i}} < 1$
and in particular that $0.5 \le \abs{M_{i} / \alpha_{i}} < 1$.
Similarly, all entries in the $j$th column of $B$ share the scale factor
\begin{equation}
  \label{eq:beta-def}
  \beta_{j} = 2^{\floor{\log_{2} N_{j}} + 1},\qquad
  N_{j}= \max_{1\le i \le k} \abs{b_{ij}},\qquad
  1 \le j \le n.
\end{equation}
Leveraging roundoff errors~\cite{ooor13}, the scaling factors in~\eqref{eq:alpha-def} and~\eqref{eq:beta-def} can be computed more efficiently in floating-point arithmetic as
\begin{equation*}
  \alpha_{i} = u^{-1} \cdot M_{i} + (1 - u^{-1}) \cdot M_{i},\qquad
  \beta_{j} = u^{-1} \cdot N_{j} + (1 - u^{-1}) \cdot M_{j}.
\end{equation*}

The scaling factors can be computed even more efficiently by relying on the bit-level representation of the floating-point numbers.
To compute $\alpha_{i}$, for example, one can take the bit string that represents $M_{i}$, set all bits that do not belong to the exponent field to zero, and add $2^{p}$ to the result using integer arithmetic.
Seen as a bit string, $2^{p}$ has only a one in position $p-1$, which corresponds to the least significant bit of the exponent field.
Therefore, adding $2^{p}$ increments the exponent by one.

For the example vectors in~\eqref{eq:ex-mats}, the scale factors are $\alpha_{1} = 2^{4}$ and $\beta_{1} = 2^{3}$, and the full block fixed-point representations are given in~\cref{fig:ex-splitting}.

\newcommand{\addz}[1]{\textcolor{Gray}{#1}}
\newcommand{\myrefh}{\rlap{$\phantom{A^{T}}$}}

\begin{figure}[t]
  \centering
  \scalebox{0.89}{%
  \begin{tikzpicture}
    [node distance=.2cm,
      every node/.style={inner sep=1pt},
      font=\footnotesize,
      every matrix/.style={inner sep=0pt, outer sep=0pt,
        left delimiter={[}, right delimiter={]},
        row sep=1pt},
      label distance=3pt]
    \matrix (AT) [matrix of math nodes,
             label=below:$A^{T}$]
      {
        \phantom{-}2^{0} \cdot 1.1001000\vphantom{\underline{0}}\\
        \phantom{-}2^{3} \cdot 1.0000000\vphantom{\underline{0}}\\
        -2^{1} \cdot 1.1101100\vphantom{\underline{0}}\\
      };
    \node [right=of AT] (tostepA1) {\myrefh$\Rightarrow 2^{4} \cdot $};
    \matrix (stepA1) [right=of tostepA1, matrix of math nodes,
             label=below:Block fixed-point]
      {
      \phantom{-}\addz{\stkout{0}.}\underline{\addz{000}}\;\underline{110}\;\underline{010}\;\underline{00\addz{0}}\vphantom{2^0}\\
      \phantom{-}\addz{\stkout{0}.}\underline{100}\;\underline{000}\;\underline{00\addz{0}}\;\underline{\addz{000}}\vphantom{2^0}\\
      -\addz{\stkout{0}.}\underline{\addz{00}1}\;\underline{110}\;\underline{110}\;\underline{0\addz{00}}\vphantom{2^0}\\
      };
    \node [right=of stepA1] (tostepA2) {\myrefh$\Rightarrow 2^{1} \cdot$};
    \matrix (stepA21) [right=of tostepA2, matrix of math nodes,
             label=below:$A^{T}_{(1)}$]
      {
      \phantom{-}000\vphantom{\underline{2}^0}\\
      \phantom{-}100\vphantom{\underline{2}^0}\\
      -001\vphantom{\underline{2}^0}\\
      };
    \node [right=of stepA21] (int21) {\myrefh$+ 2^{-2} \cdot$};
    \matrix (stepA22) [right=of int21, matrix of math nodes,
             label=below:$A^{T}_{(2)}$]
      {
      \phantom{-}110\vphantom{\underline{2}^0}\\
      \phantom{-}000\vphantom{\underline{2}^0}\\
      -110\vphantom{\underline{2}^0}\\
      };
    \node [right=of stepA22] (int22) {\myrefh$+ 2^{-5} \cdot$};
    \matrix (stepA23) [right=of int22, matrix of math nodes,
             label=below:$A^{T}_{(3)}$]
      {
      \hphantom{-}010\vphantom{\underline{2}^0}\\
      \hphantom{-}000\vphantom{\underline{2}^0}\\
      -110\vphantom{\underline{2}^0}\\
      };
    \node [right=of stepA23] (int23) {\myrefh$+ 2^{-8} \cdot$};
    \matrix (stepA24) [right=of int23, matrix of math nodes,
             label=below:$A^{T}_{(4)}$]
      {
      \phantom{-}000\vphantom{\underline{2}^0}\\
      \phantom{-}000\vphantom{\underline{2}^0}\\
      \phantom{-}000\vphantom{\underline{2}^0}\\
      };
    \matrix (BT) [matrix of math nodes,
             below=1cm of AT,
             label=below:$B$]
      {
        \phantom{-}2^{0} \cdot 1.0110001\vphantom{\underline{0}}\\
        -2^{2} \cdot 1.1110100\vphantom{\underline{0}}\\
        \phantom{-}2^{1} \cdot 1.1101000\vphantom{\underline{0}}\\
      };
    \node [right=of BT] (tostepB1) {\myrefh$\Rightarrow 2^{3} \cdot $};
    \matrix (stepB1) [right=of tostepB1, matrix of math nodes,
             label=below:Block fixed-point]
      {
      \phantom{-}\addz{\stkout{0}.}\underline{\addz{00}1}\;\underline{011}\;\underline{000}\;\underline{1\addz{00}}\vphantom{2^0}\\
      -\addz{\stkout{0}.}\underline{111}\;\underline{101}\;\underline{00\addz{0}}\;\underline{\addz{000}}\vphantom{2^0}\\
      \phantom{-}\addz{\stkout{0}.}\underline{\addz{0}11}\;\underline{101}\;\underline{000}\;\underline{\addz{000}}\vphantom{2^0}\\
      };
    \node [right=of stepB1] (tostepB2) {\myrefh$\Rightarrow 2^{0} \cdot$};
    \matrix (stepB21) [right=of tostepB2, matrix of math nodes,
             label=below:$B^{(1)}$]
      {
      \phantom{-}001\vphantom{\underline{2}^0}\\
      -111\vphantom{\underline{2}^0}\\
      \phantom{-}011\vphantom{\underline{2}^0}\\
      };
    \node [right=of stepB21] (int21) {\myrefh$+ 2^{-3} \cdot$};
    \matrix (stepB22) [right=of int21, matrix of math nodes,
             label=below:$B^{(2)}$]
      {
      \phantom{-}011\vphantom{\underline{2}^0}\\
      -101\vphantom{\underline{2}^0}\\
      \phantom{-}101\vphantom{\underline{2}^0}\\
      };
    \node [right=of stepB22] (int22) {\myrefh$+ 2^{-6} \cdot$};
    \matrix (stepB23) [right=of int22, matrix of math nodes,
             label=below:$B^{(3)}$]
      {
      \phantom{-}000\vphantom{\underline{2}^0}\\
      \phantom{-}000\vphantom{\underline{2}^0}\\
      \phantom{-}000\vphantom{\underline{2}^0}\\
      };
    \node [right=of stepB23] (int23) {\myrefh$+ 2^{-9} \cdot$};
    \matrix (stepB24) [right=of int23, matrix of math nodes,
             label=below:$B^{(4)}$]
      {
      \phantom{-}100\vphantom{\underline{2}^0}\\
      \phantom{-}000\vphantom{\underline{2}^0}\\
      \phantom{-}000\vphantom{\underline{2}^0}\\
      };
    \end{tikzpicture}}
  \caption{Bit splitting to obtain the slices for the two matrices in~\eqref{eq:ex-mats}.
On the left, the matrix entries are represented using a radix-2 scientific notation with $p=8$.
The second step uses a block fixed-point representation with a common scale and 12 significant bits.
The bits that are prepended or appended, compared with the previous step, are greyed out, and the leading bit (stricken out) is always zero.
The final step contains the slicing of each matrix into matrices with elements in $\intf[3]$.
}
\label{fig:ex-splitting}

\end{figure}

Next, we need to split these fixed-point representations into \emph{slices}.
For any $s \in \N$, the slices of $A$ and $B$ can be defined as
\begin{equation}
  \label{eq:bit-splitting-alg}
  \begin{aligned}
    A_{(\ell)} &= \left[2^{\ell t} \biggl( \diag(\alpha^{-1}) A - \sum_{r=1}^{\ell-1}2^{-rt}A_{(r)}\biggr)\right] \in \Imk,\qquad \ell = 1, 2, \ldots, s,\\
    B^{(h)} &= \left[2^{h t} \biggl( B \diag(\beta^{-1}) - \sum_{r=1}^{h-1}2^{-rt}B^{(r)}\biggr) \right] \in \Ikn, \qquad h = 1, 2, \ldots, s,
  \end{aligned}
\end{equation}
where $[\;\cdot\;]$ is the integer part operator, defined for $x \in \R$ by
\begin{equation*}
[x] =
\begin{dcases}
  \floor{x},\quad&x \ge 0,\\
  \ceil{x},  &x < 0.
\end{dcases}
\end{equation*}

We can also give a bit-level view of~\eqref{eq:bit-splitting-alg}: if we number the bits of the fraction in the block fixed-point representation from left to right, starting with index 1 to the right of the binary point, then the $\ell$th slice of $A$ contains the bits in position $(\ell-1) t+1$ to $\ell t$, and the $h$th slice of $B$ contains those in position $(h-1)t+1$ to $h t$.
In fact, this is the most natural way of understanding this slicing technique, which can be achieved by relying on bit-level operations only (bit masking and shifts).

The slicing of the example vectors in~\eqref{eq:ex-mats} is given in~\cref{fig:ex-splitting}.
In this example, setting $s = 4$ is sufficient to ensure that all bits of $A$ and $B$ are retained in split form.
The slices $A_{(4)}$ and $B^{(3)}$ contain only zeros and could in principle be ignored.
One could also choose a smaller value of $s$ and discard slices with a higher index.

\newcommand{\wt}{\widetilde}

This splitting technique yields the approximations
\begin{equation}\label{eq:slices-eq}
  \wt A = \diag(\alpha) \sum_{\ell = 1}^{s}2^{-\ell t}A_{(\ell)},\qquad
  \wt B = \sum_{h = 1}^{s}2^{-h t}B^{(h)}\diag(\beta).
\end{equation}
Using the properties of the Hadamard product, which we denote by $\circ$, one can write
\begin{equation}
  \label{eq:wtc-expr}
  \begin{aligned}
    \wt C &=  \wt A \wt B
    = \Biggl(\diag(\alpha) \sum_{\ell = 1}^{s}2^{-\ell t}A_{(\ell)}\Biggr)
       \Biggl(\sum_{h = 1}^{s}2^{-h t}B^{(h)} \diag(\beta)\Biggr)\\
    &= \alpha \beta^{T} \circ
      \sum_{\ell=1}^{s}\sum_{h=1}^{s} 2^{-(\ell+h)t}A_{(\ell)}B^{(h)}.
  \end{aligned}
\end{equation}
It is important to stress that, since $A_{(\ell)}$ and $B^{(h)}$ are both matrices of integers, the product can be computed exactly, as long as the computation and accumulation of the products is done in a wide enough format, which we call the \emph{accumulation format}.
At this point, the entries of each matrix $A_{(\ell)}B^{(h)}$ are converted
to $\fpf$,
and accumulated in floating-point arithmetic.

\newcommand{\myscale}{1.15}
\newcommand{\drawmatrix}[3]{
  \foreach \i in {1,...,#1} {
    \node at (0, \myscale*#1-\myscale*\i+\myscale) {$A_{(\i)}$};
  }
  \foreach \j in {1,...,#2} {
    \node at (\myscale*\j, \myscale*#1+0.9) {$B^{(\j)}$};
  }
  \pgfmathparse{#1 > #2 ? #1 : #2}
  \let\maxdim\pgfmathresult
  \foreach \i in {1,...,\maxdim} {
    \draw [ultra thick] (\myscale*\i,\myscale*#1) -- (\myscale,\myscale*#1-\i*\myscale+\myscale);
  }
  \foreach \j in {1,...,#2} {
    \foreach \i [evaluate={\sum=int(\i+\j);}] in {1,...,#1} {
      \ifthenelse{\sum < 6}{
        \node [fill=White,square] at (\j*\myscale,#1*\myscale-\i*\myscale+\myscale) {$2^{-\sum t}$};
      }{
        \node [dashed,fill=Gray!25,square] at (\j*\myscale,#1*\myscale-\i*\myscale+\myscale) {$2^{-\sum t}$};
      }
    }
  }
}
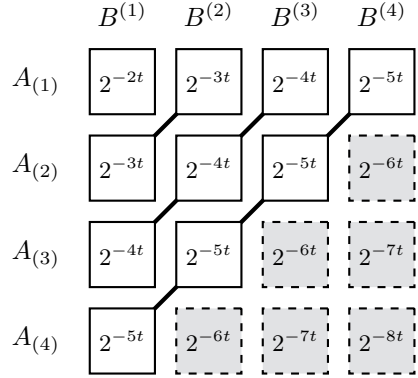
\begin{figure}[t]
  \centering
  \begin{tikzpicture}[thick,square/.style={draw,
        inner sep=0pt,font=\small,
        regular polygon,regular polygon sides=4,}]
      \drawmatrix{4}{4}{1.15}
  \end{tikzpicture}
  \caption{Products computed by different variants of the integer Ozaki scheme.
  The constant in each box is the scaling factor to be applied to the product of the slice of $A$ in the corresponding row and the slice of $B$ in the corresponding column.
  The algorithm of Ootomo, Ozaki, and Yokota~\cite{ooy24} only computes the products corresponding to boxes with a solid edge, and it accumulates them in floating-point arithmetic.
  Uchino, Ozaki, and Imamura~\cite{uoi25} use integer arithmetic to accumulate the matrices with the same scale factor (along the black diagonals) followed by accumulation of partial sums in floating-point arithmetic.}
  \label{fig:products}
\end{figure}

Note that in~\eqref{eq:wtc-expr}, the contribution of the term $A_{(\ell)}B^{(h)}$ to the final result is inversely proportional to the sum $\ell + h$, since the multiplier $2^{-(\ell+h)t}$ will be small.
For this reason, Ootomo, Ozaki, and Yokota~\cite{ooy24} have proposed to only compute $A_{(\ell)}B^{(h)}$ if $\ell + h \le s+1$.
This is similar to the strategy used for matrix multiplication in multi-word arithmetic~\cite{fhlm23}.
The matrix products computed by this algorithm are shown in \cref{fig:products}.

How many bits should the accumulation format have?
To store exactly the product of two $(t+1)$-bit integers, $2t+1$ bits are necessary, and in view of~\eqref{eq:max-num-bits-sum-n} we need at most $\ceil{\log_{2} k}$ additional bits to store exactly the sum of the $k$ partial products.
Therefore, the sum can be represented exactly as long as the output format $\intf[T]$ has
\begin{equation*}
  T = (2t+1) + \ceil{\log_{2}k} - 1 = 2t + \ceil{\log_{2}k}
\end{equation*}
bits, where $-1$ is needed because we assume that $\intf[T]$ is a signed integer format.
For the example in~\cref{fig:ex-splitting}, $t = 3$ and $k = 3$ imply that the accumulation format requires $T > 8$.
The alignment of the partial products in the final accumulation is shown in \cref{fig:bit-alignment}.

\newlength{\myhsep}
\setlength{\myhsep}{11pt}
\newcommand{\mysep}[3][solid]{%
    \draw[#1] ($(#2)!0.5!(#3)+(-1,0)$) -- ($(#2)!0.5!(#3)+(8,0)$);}

\begin{figure}[t]
  \centering
  \begin{tikzpicture}[node distance=.1cm,
      every node/.style={minimum width=1.5cm,inner sep=1pt}]
    \node (P11) {$A_{(1)}B^{(1)}$};
    \node[right=\myhsep of P11] {$2^{1\hphantom{0-}} \cdot\hphantom{--} \llap{$-$}00011111$}; %

    \node[below=of P11] (P12) {$A_{(1)}B^{(2)}$};
    \node[below=of P12] (P21) {$A_{(2)}B^{(1)}$};
    \node[right=\myhsep of P12] {$2^{-2\hphantom{0}} \cdot\hphantom{--} \hphantom{000}\llap{$-$}00011001$}; %
    \node[right=\myhsep of P21] {$2^{-2\hphantom{0}} \cdot\hphantom{--} \hphantom{000}\llap{$-$}00001100$}; %

    \node[below=of P21] (P13) {$A_{(1)}B^{(3)}$};
    \node[below=of P13] (P22) {$A_{(2)}B^{(2)}$};
    \node[below=of P22] (P31) {$A_{(3)}B^{(1)}$};
    \node[right=\myhsep of P13] {$2^{-5\hphantom{0}} \cdot\hphantom{--} \hphantom{000000}00000000$}; %
    \node[right=\myhsep of P22] {$2^{-5\hphantom{0}} \cdot\hphantom{--} \hphantom{000000}\llap{$-$}00001100$}; %
    \node[right=\myhsep of P31] {$2^{-5\hphantom{0}} \cdot\hphantom{--} \hphantom{000000}\llap{$-$}00010000$}; %

    \node[below=of P31] (P14) {$A_{(1)}B^{(4)}$};
    \node[below=of P14] (P23) {$A_{(2)}B^{(3)}$};
    \node[below=of P23] (P32) {$A_{(3)}B^{(2)}$};
    \node[below=of P32] (P41) {$A_{(4)}B^{(1)}$};
    \node[right=\myhsep of P14] {$2^{-8\hphantom{0}} \cdot\hphantom{--} \hphantom{000000000}00000000$}; %
    \node[right=\myhsep of P23] {$2^{-8\hphantom{0}} \cdot\hphantom{--} \hphantom{000000000}00000000$}; %
    \node[right=\myhsep of P32] {$2^{-8\hphantom{0}} \cdot\hphantom{--} \hphantom{000000000}\llap{$-$}00011000$}; %
    \node[right=\myhsep of P41] {$2^{-8\hphantom{0}} \cdot\hphantom{--} \hphantom{000000000}00000000$}; %

    \node[below=of P41] (P24) {$A_{(2)}B^{(4)}$};
    \node[below=of P24] (P33) {$A_{(3)}B^{(3)}$};
    \node[below=of P33] (P42) {$A_{(4)}B^{(2)}$};
    \node[right=\myhsep of P24] {$2^{-11} \cdot\hphantom{--} \hphantom{000000000000}00011000$}; %
    \node[right=\myhsep of P33] {$2^{-11} \cdot\hphantom{--} \hphantom{000000000000}00000000$}; %

    \node[right=\myhsep of P42] {$2^{-11} \cdot\hphantom{--} \hphantom{000000000000}00000000$}; %

    \node[below=of P42] (P34) {$A_{(3)}B^{(4)}$};
    \node[below=of P34] (P43) {$A_{(4)}B^{(3)}$};
    \node[right=\myhsep of P34] {$2^{-14} \cdot\hphantom{--} \hphantom{000000000000000}00001000$}; %
    \node[right=\myhsep of P43] {$2^{-14} \cdot\hphantom{--} \hphantom{000000000000000}00000000$}; %

    \node[below=of P43] (P44) {$A_{(4)}B^{(4)}$};
    \node[right=\myhsep of P44] {$2^{-17} \cdot\hphantom{--} \hphantom{000000000000000000}00000000$}; %

    \node[below=of P44] (result) {$AB_{\vphantom{(4)}}^{\vphantom{(4)}}$};
    \node[right=\myhsep of result] {$2^{-17} \cdot\hphantom{--} \llap{$-$}00100100000110100111000000$};

    \mysep[dashed]{P11}{P12}
    \mysep[dashed]{P21}{P13}
    \mysep[dashed]{P31}{P14}
    \mysep[solid]{P41}{P24}
    \mysep[dashed]{P42}{P34}
    \mysep[dashed]{P43}{P44}
    \mysep[ultra thick]{P44}{result}

  \end{tikzpicture}
\caption{Alignment of bits in the 16 products of the form $A_{(\ell)}B^{(h)}$ for the slices in~\cref{fig:ex-splitting}.
The dashed lines separate blocks of partial products with the same scale factor, which lie along the same diagonal in \cref{fig:products}.
The products below the thin, solid line correspond to the greyed-out boxes with a dashed border in~\cref{fig:products}.
The value below the thick solid line is the full-precision fixed-point representation of the result including all products.
In this case, this is the exact result, because all the bits in $A$ and $B$ were allocated to a slice, and all slices were used in the computation.}
\label{fig:bit-alignment}
\end{figure}

\begin{figure}[t]
  \centering
  \begin{tikzpicture}[scale=0.5]

    \draw[draw] (0,0) rectangle (20,1);
    \draw[draw, fill=Gray!20] (0,0) rectangle (8,1);
    \draw[draw, fill=Cyan!20] (8,0) rectangle (12,1);
    \draw[draw, fill=BrickRed!20] (12,0) rectangle (20,1);

    \node at  (4, -0.75)  {\vphantom{d$\ceil{\log_{2}k}$}unused};
    \node at  (10, -0.75)  {\vphantom{d$\ceil{\log_{2}k}$}$\ceil{\log_{2}k}$};
    \node at (16, -0.75) {\vphantom{d$\ceil{\log_{2}k}$}$2t+1$};

    \draw [decorate,decoration={brace,amplitude=10pt}] (0,1.2) -- (20,1.2) node [midway, above=10pt] {$T+1\;$};
    \draw [decorate,decoration={brace,amplitude=10pt}] (20,-1.2) -- (8,-1.2) node [midway, below=10pt] {$T'+1\;$};
    \draw [decorate,decoration={brace,amplitude=10pt}] (12,-2) -- (0,-2) node [midway, below=10pt] {$K$};

  \end{tikzpicture}
\caption{Allocation of bits in the accumulation format $\intf[T]$ for the sum of $k$ products of values in $\intf$.
  The integers $K$ and $T'$ are defined in~\eqref{eq:optimal-K} and~\eqref{eq:used-acc-bits}, respectively.}
  \label{fig:bit-allocation}
\end{figure}
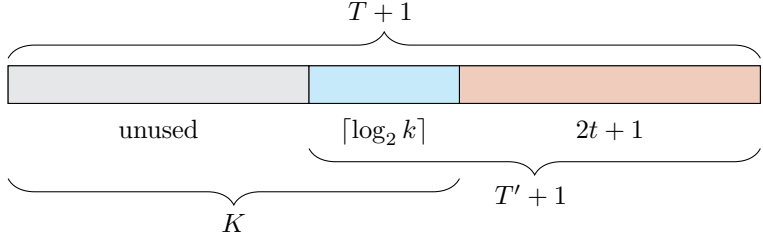

In practice, however, the width of the accumulation format depends on what is already present in the hardware.
It is therefore more appropriate to assume that the hardware MMA unit accepts inputs in $\intf[t']$ and accumulates and returns outputs in $\intf[T]$, and ask what is the maximum number of bits per slice.
If the input to the MMA has at most $t'+1$ bits in two's complement, then we must have that $t \le t'$, but to ensure that the sum of $k$ products of integers in $\intf$ can be represented in $\intf[T]$, we must also require that
\begin{equation*}
    t \le \biggl\lfloor\frac{T - \ceil{\log_{2}k}}{2} \biggr\rfloor.
\end{equation*}
Therefore, the optimal choice for $t$ is
\begin{equation}
  \label{eq:num-int-bit}
   t = \min\left\{t', \biggl\lfloor\frac{T - \ceil{\log_{2}k}}{2} \biggr\rfloor\right\},
\end{equation}
as this value maximizes the number of bits per slice, and therefore reduces the overall number of slices needed to satisfy a given accuracy threshold.
Yet another way to look at this question is to ask what is the largest $k$ that is allowed by the algorithm assuming $\intf[t']$ and $\intf[T]$ as input and output format of the integer MMA.
In view of the discussion above, this will be
\begin{equation}
  \label{eq:optimal-K}
  k = 2^{K},\qquad K = T + 1 - 2(t'+1) = T - 2t' + 1.
\end{equation}
Here $K$ represents the number of bits that can be used for the accumulation once the bits for a single product have been accounted for.
\Cref{fig:bit-allocation} shows how the bits in the accumulation format are allocated.
As $t$ cannot be smaller than 1, this algorithm will only work provided that $k \le 2^{T - 3}$.
\Cref{alg:ooy-24} summarizes the method.

One example of existing hardware suitable for the integer Ozaki scheme are the integer tensor cores, which %
use $\intf[31]$ (INT32) as accumulation format, and either $\intf[7]$ (INT8) or $\intf[3]$ (INT4) for the inputs.
According to \eqref{eq:optimal-K}, with these combinations of input and output formats, the largest value of $k$ the algorithm can support is $2^{16}=65{,}536$ for~INT8 and $2^{24}=16{,}777{,}216$ for~INT4.

\begin{algorithm2e}[t]
  \caption{Matrix multiplication using an integer MMA unit~\cite{ooy24}.}
  \label{alg:ooy-24}
  \KwIn{$A \in \Fmk$, $B \in \Fkn$, $s \in \N$, MMA with $\intf[t']$ inputs and $\intf[T]$ outputs.}
  \KwOut{$C \in \Fmn$ such that $C \approx AB$}
   $t \gets \min\left\{t',\big\lfloor(T - \ceil{\log_{2}k})/2\big\rfloor\right\}$\; \label{ln:compute-t}
  \For{$i \gets 1\ \textbf{to}\ m$}{
    $M \gets \max_{1 \le j \le k}\abs{a_{ij}}$\;
    $\alpha_{i} \gets 2^{\floor{\log_{2}M}+1}$\;
  }
  \For{$\ell \gets 1\ \textbf{to}\ s$}{
    $\mathcal{A} \gets $ block fixed-point representation of $A$ with row $i$ scaled by $\alpha_{i}$.\;
    $A_{(\ell)} \gets $ bits from $(\ell - 1) t +1$ to $\ell t$ of $\mathcal{A}$.
  }
  \For{$j \gets 1\ \textbf{to}\ n$}{
    $N \gets \max_{1 \le i \le k}\abs{b_{ij}}$\;
    $\beta_{i} \gets 2^{\floor{\log_{2}N}+1}$\;
  }
  \For{$h \gets 1\ \textbf{to}\ s$}{
    $\mathcal{B} \gets $ block fixed-point representation of $B$ with column $i$ scaled by $\beta_{i}$.\;
    $B^{(h)} \gets $ bits from $(h - 1)t +1$ to $h t$ of $\mathcal{B}$.
  }
  $S \gets 0_{m \times n}$\;
  \For{$\ell \gets 1\ \textbf{to}\ s$}{
    \For{$h \gets 1\ \textbf{to}\ s-\ell+1$}{
      $E \gets A_{(\ell)}B^{(h)} \text{converted to } \fpf$.\;
      $S \gets S + 2^{-(\ell+h)t}E$\; \label{ln:accummulation}
    }
  }
  $ C \gets \diag(\alpha) S \diag(\beta)$
\end{algorithm2e}

\paragraph{Enhanced algorithm}

More recently, Uchino, Ozaki, and Imamura~\cite{uoi25} have proposed two improvements on the integer Ozaki scheme just described.

The first tweak is aimed at improving the performance of the algorithm on hardware equipped with fast integer MMA units.
By profiling the original integer implementation~\cite{ooy24}, the authors discover that a significant portion of the runtime is spent accumulating the integer matrix products in floating-point arithmetic.
This operation is slow because, unlike the integer matrix products, it cannot be performed by the efficient MMA units and relies on the general purpose floating-point units on the GPU.
Therefore, they suggest accumulating the matrix products on line~\ref{ln:accummulation} of \cref{alg:ooy-24} in integer arithmetic for all pairs of indices $(\ell, h)$ such that $\ell + h$ is constant.
This strategy, depicted in \cref{fig:products} reduces the number of floating-point sums per element of the result from $s(s+1) / 2$ (the number of gray boxes in the figure) to $s-1$ (the number of diagonals).
This will increase the number of integer sums from $k$ to $k + s - 1$, and if the parameter $t$ is computed as in \cref{ln:compute-t}, this may cause some of the integer sums to overflow.
The simplest way to address this would be to change the definition of $t$ to
\begin{equation*}
    t \gets \min\left\{t',\biggl\lfloor\frac{T - \ceil{\log_{2}(k + s - 1)}}{2}\biggr\rfloor\right\}.
\end{equation*}
Uchino, Ozaki, and Imamura suggest a more refined technique that reverts to floating-point arithmetic when the maximum number of error-free sums for the choice of $t$ in \cref{ln:compute-t} is reached.
This is a more effective solution in practice, since reducing $t$ is very likely to increase the number of splits required, while floating-point accumulation will only be necessary when $k \lesssim 2^{T - 2t - \lceil \log_2 k \rceil}$, in which case $k \gg s$.

The second enhancement pertains to the technique used to slice the fixed-point representation into integers.
The bit splitting technique in~\cref{fig:ex-splitting} implicitly uses round-to-zero, as shown by~\eqref{eq:bit-splitting-alg}.
When the fixed-point representation is truncated by fixing the number of slices, using round-to-nearest may yield a more accurate result.
This splitting technique is more expensive to implement than truncation and requires that the first slice have $t-1$ rather than $t$ bits.

\section{Error analysis}
\label{sec:error-analysis}

The Ozaki scheme is subject to two sources of error.
A \emph{truncation error} is incurred by approximating $A$ and $B$ with $\wt A$ and $\wt B$, respectively.
The magnitude of this error depends on the number $s$ of slices in~\eqref{eq:slices-eq}.
The second source are \emph{rounding errors}, due to the use floating-point arithmetic to accumulate the intermediate products of the form $A_{(\ell)}B^{(h)}$, which are computed exactly unless overflow occurs.

These two types of error arise in different parts of the computation: truncation errors are caused by the slicing of $A$ and $B$, while rounding errors are caused by the use of floating-point arithmetic in the final accumulation.
Therefore, we treat these two sources of error separately.
\Cref{sec:ea-truncation-error,sec:ea-all-products,sec:ea-reduced-products} deal with truncation errors, while rounding errors are the main subject of \cref{sec:ea-fp-accumulation}.
Finally \cref{sec:ea-discussion} combines the two into new error bounds on matrix multiplication.

We remark that Uchino, Ozaki, and Imamura~\cite{uoi25} have already undertaken error analysis of these algorithms.
Our discussion is different in several respects.
\begin{itemize}
  \item It gives precise conditions under which the algorithm may fail even when a large number of slices are used.
 This theoretical result is confirmed numerically by the experiments in \cref{sec:minimal-example}.

  \item It allows for a different number of slices for $A$ and $B$, thereby enabling the development of more flexible algorithms.
    In particular, we show that if one of the two matrices is \emph{badly scaled} but the other is not, then reducing the number of slices used for the well scaled matrix does not affect the accuracy of the result.
    This is confirmed by the numerical experiments in \cref{sec:lu-experiment}.

  \item It yields data-dependent bounds that can be used in practice to decide how many slices should be used for a specific choice of $A$ and $B$.

  \item It improves the observation in~\cite[sect.~5.1]{uoi25} that depending on the floating-point format of input and output matrices and on the integer formats used by the MMA unit, some of the sums in floating-point arithmetic can be computed exactly.
  Specifically, we explain in which order the partial products can be accumulated to maximize accuracy, and we provide a bound for the maximum error of this strategy.

  \item It relies on error analysis in the style of Wilkinson~\cite{wilk63} and Higham~\cite{high02}, whilst the previous work used the error bounds of Jeannerod and Rump~\cite{jeru13}.
\end{itemize}

\subsection{General results on truncation error}
\label{sec:ea-truncation-error}

We start by discussing the error in the splitting~\eqref{eq:slices-eq}.
Unlike the existing algorithms~\cite{ooy24},~\cite{uoi25}, we allow for different numbers of splits in the two input matrices.
Since the entries of $A$ and $B$ are floating-point numbers with finite precision and limited exponent range, the number of nonzero slices must be finite.
Let $\sAs, \sBs \in \N$ be the smallest integers that satisfy
\begin{equation}
  \label{eq:exact-splitting}
  A = \diag (\alpha) \sum_{\ell=1}^{\sAs} 2^{-\ell t} A_{(\ell)}\quad\text{and}\quad
  B = \sum_{h=1}^{\sBs} 2^{-h t} B^{(h)} \diag (\beta).
\end{equation}
In other words, $\sAs$ and $\sBs$ are the minimum number of slices required to represent exactly the smallest element of $A$ and $B$, respectively, in the fixed-point representation in \cref{fig:ex-splitting}.
Using $\sAs - 1$ or fewer slices would not satisfy the equality on the left of~\eqref{eq:exact-splitting}, and using $\sBs - 1$ or fewer slices would invalidate the equality on the right.

The cost of the algorithm depends directly on the number of slices of $A$ and $B$ used, and $\sAs$ and $\sBs$ may, in general, be too large for the algorithm to be practical.
Therefore, we consider what happens when $A$ and $B$ are split into $\sA < \sAs$ and $\sB < \sBs$ slices, respectively.

The $i$th row of $A$ is converted to a block with scale $\alpha_{i}$ and $\sA t$ bits to the right of the binary point.
Therefore, the absolute conversion error of entry $a_{ij}$ must be bounded, in magnitude, by $2^{-\sA t}\alpha_{i}$, and we can write
\begin{equation}\label{eq:a-split}
    A = \DA + \diag(\alpha) \sum_{\ell = 1}^{\sA}2^{-\ell t}A_{(\ell)},\qquad
    \abs{\da_{ij}} < \alpha_{i} 2^{-\sA t},
\end{equation}
where the entries of $\alpha \in \N[m]$ are defined in~\eqref{eq:alpha-def}, and $A_{(\ell)} \in \intf[{\mbyk}]$ for $\ell = 1, 2, \ldots, \sA$.
Similarly, for the $j$th column of $B$ we have
\begin{equation}\label{eq:b-split}
  B = \DB + \sum_{h = 1}^{\sB}2^{-h t}B^{(h)} \diag(\beta),\qquad
  \abs{\db_{ij}} < \beta_{j}2^{-\sB t},
\end{equation}
where the entries $\beta \in \N[n]$ are defined in~\eqref{eq:beta-def} and $B^{(h)} \in \intf[{\kbyn}]$ for $h = 1, 2, \ldots, \sB$.

Equations~\eqref{eq:a-split} and~\eqref{eq:b-split} bound the absolute conversion error, but, in error analysis, it is often more informative to bound the relative error instead.
Such bounds arise naturally when using floating-point arithmetic, because floating-point numbers have constant precision.
In fixed-point arithmetic, smaller numbers have lower precision, and bounds like the one in~\eqref{eq:a-split} and~\eqref{eq:b-split} are therefore more familiar.

To obtain a relative bound, we note that
\begin{align}
  \frac{\abs{\da_{ij}}}{\abs{a_{ij}}}
  &< \frac{\alpha_{i}}{\abs{a_{ij}}} 2^{-st} \label{eq:a-bound-element}\\
  &\le \frac{2 \max_{j} \abs{a_{ij}}}{\min_{j}\abs{a_{ij}}} 2^{-st} \label{eq:a-bound-row}\\
  &\le \kA 2^{-st},
  &\kA := 2 \max_{i}\frac{\max_{j} \abs{a_{ij}}}{\min_{j}\abs{a_{ij}}}
    \label{eq:a-bound-row-2}.
\end{align}
The bound~\eqref{eq:a-bound-element} can be large if $a_{ij}$ is small in magnitude.
The bound~\eqref{eq:a-bound-row} is saying that if the $i$th row of $A$ is \emph{badly scaled}, that is, has entries that vary widely in magnitude, then the conversion error can be large, in relative terms, for the entries of that row.
In fact, for a badly scaled row, the quantity
\begin{equation*}
  \frac{\max_{j} \abs{a_{ij}}}{\min_{j}\abs{a_{ij}}}
\end{equation*}
will necessarily be large, and the fact that
\begin{equation*}
  \max_{j} \abs{a_{ij}} < \alpha_{i} \le 2 \max_{j} \abs{a_{ij}}
\end{equation*}
shows that the relative error itself is likely to be large for some entries of that row.
Similarly, the bound~\eqref{eq:a-bound-row-2} will only be large if the matrix has at least one badly scaled row, whose entries are prone to a large conversion error in relative terms.

For the conversion error in~\eqref{eq:b-split}, we can look at the columns of $B$ to obtain the bound
\begin{align}
  \label{eq:b-bound}
  \frac{\abs{\db_{ij}}}{\abs{b_{ij}}} < \frac{\beta_{j}}{\abs{b_{ij}}} 2^{-st} \le
  \frac{2 \max_{i} \abs{b_{ij}}}{\min_{i}\abs{b_{ij}}} 2^{-st} \le \kB 2^{-st}, \qquad
  \kB := 2 \max_{j}\frac{\max_{i} \abs{b_{ij}}}{\min_{i}\abs{b_{ij}}}.
\end{align}

\subsection{Computation of all products}
\label{sec:ea-all-products}

We begin by considering a variant of \cref{alg:ooy-24} that computes all $\sA \sB$ products without the early termination demonstrated in \cref{fig:products}.
Using~\eqref{eq:a-split} with $s = \sA$ and~\eqref{eq:b-split} with $s = \sB$, the matrix $C:=AB$ can be written as
\begin{equation}
  \label{eq:c-expr}
  \begin{aligned}
    C
    &= \biggl(\DA + \diag(\alpha) \sum_{\ell = 1}^{\sA}2^{-\ell t}A_{(\ell)}\biggr)
      \biggl(\DB + \sum_{h = 1}^{\sB}2^{-h t}B^{(h)} \diag(\beta)\biggr)\\
    &= \DA \cdot B + A \cdot \DB + \DA \cdot \DB + \wt C,
  \end{aligned}
\end{equation}
for $\wt C$ in~\eqref{eq:wtc-expr}.
By construction, the entries of each partial product $A_{(\ell)}B^{(h)}$ are stored in an integer format $\intf[T]$
that is sufficiently large to guarantee that the matrix products are exact.
These $\sA \sB$ integer matrices are accumulated in floating-point arithmetic, and they must therefore be converted to a floating-point format, $\fpf[p]$ say, before the accumulation can take place.
If $p \ge T$, then the significand of the floating-point numbers are large enough to store the partial results to full precision, and the conversion from $\intf[T]$ to $\fpf[p]$ will be exact.
The following scaling is also exact, since the entries of $\alpha$ and $\beta$ are all powers of two.

Let $\wh C$ be the matrix obtained by accumulating the integer matrices in floating-point arithmetic.
Standard error analysis of floating-point summation~\cite[Chap.~3]{high02} gives the bound
\begin{equation}\label{eq:c-hat}
  \wh C
  = (1+ \Theta_{ij}) \circ \wt C,
  \qquad \abs{\theta_{ij}} \le \gamma_{\psi - 1}
\end{equation}
where
\begin{equation}
  \label{eq:gamma-def}
  \gamma_{n} := \frac{nu}{1-nu},
\end{equation}
$\psi$ is the number of matrices to be added, and $u := 2^{-p}$ is the unit roundoff of $\fpf[p]$.
The bound~\eqref{eq:c-hat} certainly holds for $\psi = \sA \sB$, but, as noted in~\cite[sect.~5.1]{uoi25}, some of these $\sA \sB - 1$ sums can be computed without rounding errors in floating-point arithmetic, when $T$ is sufficiently smaller than $p$.
Therefore, in this section we keep~$\psi$ generic, and refer the reader to \cref{sec:ea-fp-accumulation} for a discussion of the exact value of~$\psi$ for different variants of the algorithm.

Combining~\eqref{eq:c-expr} and~\eqref{eq:c-hat} gives
\begin{equation}\label{eq:err-bound-partial}
  \begin{aligned}
    \abs{C - \wh C}
    &\le \abs{\DA} \abs{B} + \abs{A} \abs{\DB} + \abs{\DA} \abs{\DB} +
      \gamma_{\psi - 1} \abs{\wt C}\\
    &\le \zeta_{A,B} \abs{A} \abs{B} + \gamma_{\psi - 1} \abs{\wt C},
  \end{aligned}
\end{equation}
where
\begin{equation}
  \label{eq:zeta-ab}
  \zeta_{A,B} :=  2^{-\sA t}\kA + 2^{-\sB t}\kB + 2^{-(\sA+\sB) t} \kA \kB
\end{equation}
and the absolute value of a matrix is to be understood entry-wise.
To write the whole bound in terms of $\abs{A}\abs{B}$, we can use the fact that $\wt C = (A - \DA)(B - \DB)$ to obtain
\begin{equation}
  \label{eq:d-bound}
  \abs{\wt C} \le (1 + \zeta_{A,B}) \abs{A} \abs{B}.
\end{equation}
Finally, plugging~\eqref{eq:d-bound} into~\eqref{eq:err-bound-partial} gives the error bound
\begin{equation}
  \label{eq:bound-full}
  \abs{\wh C - C}
  \le \bigl(\zeta_{A,B} + \gamma_{\psi - 1}(1 + \zeta_{A,B})
  \bigr) \abs{A} \abs{B}.
\end{equation}

What does~\eqref{eq:bound-full} tell us?
If $2^{-\sA t}$ and $2^{-\sB t}$ are smaller than the unit roundoff of $\fpf[p]$ in which the accumulation is performed, then we can write the first-order approximation of~\eqref{eq:bound-full} as
\begin{equation}
  \label{eq:first-order-bound}
  \abs{\wh C - C} \lesssim \left(
    2^{-\sA t}\kA + 2^{-\sB t}\kB + (\psi - 1) 2^{-p}
  \right) \abs{A} \abs{B}.
\end{equation}
The quantities $2^{-\sA t}$ and $2^{-\sB t}$ depend directly on the number of bits used in the block fixed-point representation of $A$ and $B$, respectively, and can therefore be taken as a measure of the limiting accuracy in the sliced matrices.
Therefore, \eqref{eq:first-order-bound} is saying that the error can potentially be large if $A$ or $B$ are badly scaled.
A large value of $\kA$ or $\kB$ can be balanced by increasing $\sA$ and $\sB$ accordingly.
In most practical scenarios, $\sA$ and $\sB$ will be small, and will therefore only moderately affect the third term in~\eqref{eq:first-order-bound}, but will have a major impact on the performance of the algorithm.

\subsection{Reduction in number of products}
\label{sec:ea-reduced-products}

We now consider the approach in \cref{alg:ooy-24}, where matrix products of the form $A_{(\ell)}B^{(h)}$ are not computed if \mbox{$\ell + h > s + 1$}.
If we extend this approach to the case $\sA \neq \sB$, then the product $A_{(\ell)}B^{(h)}$ will be computed if $\ell + h \le \max(\sA, \sB) + 1$, provided that $\ell \le \sA$ and $h \le \sB$.
We discuss in detail the case $\sA \le \sB$, but a bound for $\sA > \sB$ can be obtained in an analogous way.
If we set
\begin{equation}
  \label{eq:c-prime}
  C' = \alpha \beta^{T} \circ
  \sum_{\ell=1}^{\sA}\sum_{h=1}^{\sB - \ell + 1} 2^{-(\ell+h)t}A_{(\ell)}B^{(h)},
\end{equation}
then we can rewrite~\eqref{eq:c-expr} as
\begin{equation}
  \label{eq:red-prod}
  C = \DA \cdot B + A \cdot \DB + \DA \cdot \DB + C' + (\wt C - C').
\end{equation}
We can bound the magnitude of the last term on the right hand side of~\eqref{eq:red-prod} by noting that
\begin{equation}
  \label{eq:bound-cprime-derivation}
  \begin{aligned}
    \abs{\wt C - C'}
    &=\Bigg\lvert\diag(\alpha)
    \sum_{\ell = 1}^{\sA} \sum_{h = \sB - \ell + 2}^{\sB} 2^{-(\ell+h)t}A_{(\ell)}B^{(h)}
    \diag(\beta)\Bigg\rvert\\
    &= \Bigg\lvert\diag(\alpha)
      \sum_{\ell = 1}^{\sA} 2^{{-\ell t}} A_{(\ell)}
      \sum_{h = \sB - \ell + 2}^{\sB} 2^{-h t}B^{(h)}
      \diag(\beta) \Bigg\rvert\\
    &\le\diag(\alpha)
      \sum_{\ell = 1}^{\sA} 2^{{-\ell t}} \big\lvert A_{(\ell)} \big\rvert
      \Bigg\lvert
      \sum_{h = \sB - \ell + 2}^{\sB} 2^{-h t}B^{(h)}
      \diag(\beta) \Bigg\rvert\\
    &\le \diag(\alpha)
      \sum_{\ell = 1}^{\sA} 2^{{-\ell t}} \big\lvert A_{(\ell)} \big\rvert
      \kB 2^{-(\sB-\ell+1)t}\abs{B}\\
    &\le 2^{-\sB t} \kB \diag(\alpha)
      \sum_{\ell = 1}^{\sA} 2^{-t} \big\lvert A_{(\ell)} \big\rvert  \abs{B}\\
    &\le 2^{-\sB t} \sA \kA \kB \abs{A} \abs{B}.
  \end{aligned}
\end{equation}
After rewriting the expression in a more convenient form, the third step uses the triangular inequality, the fourth relies on~\eqref{eq:b-split} with~\eqref{eq:b-bound}, after noting that
\begin{equation*}
  \sum_{h = \sB - \ell + 2}^{\sB} 2^{-h t}B^{(h)}\diag(\beta) =
  B - \sum_{h = 1}^{\sB -\ell + 1} 2^{-h t}B^{(h)}\diag(\beta).
\end{equation*}
The last step exploits the fact that the entries of $A_{\ell}$ are bounded in magnitude by $2^{t} - 1$, combined with the observation that $\diag(\alpha) 1_{m \times k} \le \kA \abs{A}$, where $1_{m \times k}$ is the $m\times k$ matrix of ones---this is a consequence of~\eqref{eq:a-bound-row-2} for $s = 0$.

Now, let $\wh C'$ be the matrix obtained by accumulating the matrix products in~\eqref{eq:c-prime} in floating-point arithmetic.
A simple calculation shows that, in this case, the number of products to be computed is
\begin{equation*}
    \chi(\sA, \sB) = \frac{\sm(2 \sM - \sm + 1)}{2},\qquad\sM=\max(\sA,\sB),\qquad\sm=\min(\sA,\sB).
\end{equation*}
Therefore the bound
\begin{equation}
  \label{eq:c-hat-reduced}
  \wh C'
  = (1+ \Theta_{ij}) \circ C',
  \qquad \abs{\theta_{ij}} \le \gamma_{\psi - 1},
\end{equation}
where $\gamma_{n}$ is defined in~\eqref{eq:gamma-def}, would hold for $\psi = \chi(\sA,\sB)$, but some of the sums are exact in floating-point arithmetic, and we refer the reader to \cref{sec:ea-fp-accumulation} for a discussion of this.

By combining~\eqref{eq:red-prod},~\eqref{eq:bound-cprime-derivation}, and~\eqref{eq:c-hat-reduced}, we obtain the bound
\begin{equation*}
  \abs{C - \wh C'} \le (\zeta_{A,B} + 2^{-\sB t} \sA \kA \kB ) \abs{A}\abs{B}
  + \gamma_{\psi} \abs{C'},
\end{equation*}
where $\zeta_{A,B}$ is defined in~\eqref{eq:zeta-ab}.
As done in~\eqref{eq:d-bound}, we can express $\abs{C'}$ in terms of $\abs{A}$ and $\abs{B}$ by solving~\eqref{eq:red-prod} for $C'$ and taking the absolute value.
This yields the bound
\begin{equation*}
  \abs{C - \wh C'} \le \bigl(\zeta_{A,B} + 2^{-\sB t} \sA \kA \kB + \gamma_{\psi}(1 + \zeta_{A,B} + 2^{-\sB t} \sA \kA \kB)\bigr) \abs{A}\abs{B}.
\end{equation*}
This bound will grow faster than~\eqref{eq:bound-full} when $\kA$ is large, since $\zeta_{A,B}$ only features terms in $2^{-\sB t} \kB$.
When $\sA > \sB$, we can switch the role of $A$ and $B$ in \eqref{eq:bound-cprime-derivation} and obtain the equivalent bound
\begin{equation*}
  \abs{C - \wh C'} \le \bigl(\zeta_{A,B} + 2^{-\sA t} \sB \kA \kB + \gamma_{\psi}(1 + \zeta_{A,B} + 2^{-\sA t} \sB \kA \kB)\bigr) \abs{A}\abs{B},
\end{equation*}
which can grow faster than~\eqref{eq:bound-full} for large $\kappa_{B}$, as $\zeta_{A,B}$ only contains terms in $2^{-\sA t} \kA$.

\newcommand{\C}[1]{\ensuremath{C_{#1}}}
\subsection{Accumulation in floating-point arithmetic}
\label{sec:ea-fp-accumulation}

We now discuss the best strategies to accumulate the (exact) partial products of the form $A_{(\ell)}B^{(h)}$ in floating-point arithmetic.
We assume that $T < p$, where $T$ is the number of bits in the output format of the integer accumulator and $p$ is the precision of the floating-point format used.
We further define the quantity
\begin{equation}
  \label{eq:used-acc-bits}
  T' := 2t + \ceil{\log_{2}k},
\end{equation}
which represents the number of bits in the accumulator format that are actually used in the computation.

In our analysis, we divide the sums performed in the final accumulation step into several \emph{levels}.
The sum within each level will be computed exactly, and each level will produce only one matrix computed exactly.
Rounding errors will therefore only occur when adding across levels, and the number of levels will give us the constants $\psi$ to use in~\eqref{eq:c-hat} and $\psi$ to use in~\eqref{eq:c-hat-reduced}.

In order to simplify the terminology, we will refer to \cref{fig:products} and explain how the partial products can be accumulated proceeding by diagonals.
Our algorithm starts from the top left diagonal in \cref{fig:products}, which contains only the product \mbox{$\C{0} := A_{(1)}B^{(1)}$}.
The elements of this matrix belong to $\intf[T']$, and since $T' \le T \le p$, they can be converted to $\fpf$ exactly.

We then move to the second diagonal, which contains two elements, $A_{(1)}B^{(2)}$ and $A_{(2)}B^{(1)}$.
The entries of the individual products can be converted exactly to the floating-point format, and in order to compute the sum $\C{1} := A_{(1)}B^{(2)} + A_{(2)}B^{(1)}$ exactly, we only need that $T \ge T' + 1$, to allow for the possible carry.
In fact, matrix entries occupying the same position in the two products have the same exponent, and their fractions are therefore aligned.

Next, we need to add $\C{0}$ and $\C{1}$.
Floating-point addition is performed by converting the values to be added to a fixed-point format, and this is achieved by shifting right the fraction of the smaller value in magnitude, so that the two summands have the same exponent.
This alignment will shift the elements of $\C{1}$ to the right by $t$ places, and these $t$ places will be more than sufficient to store the carry bit produced when computing $\C{1}$.
We will use later the fact that, because of this shift, diagonal 1 could, in principle, accommodate another $2^{t}-2$ carries produced by subsequent diagonals.
We say that diagonal 1 has $2^{t}-2$ \emph{spare carry locations}.
The floating-point format will also need an additional bit for the carry, so that to be represented exactly the sum $\C{0} + \C{1}$ requires that $p \ge T' + t + 1$.

Accumulating the products along diagonal 2 will produce $\C{2}:= A_{(1)}B^{(3)} + A_{(2)}B^{2} + A_{(3)}B^{(1)}$, which is the sum of three partial products.
Computing this sum will require two additional carry bits, covered by the $2t$ bits by which $\C{2}$ has to be shifted.
Therefore, representing $\C{0} + \C{1} + \C{2}$ exactly requires that $p \ge T' + 2t + 1$, and we also note that diagonal 2 has $2^{t}-3$ spare carry locations.
In general, accumulating the diagonal products up to
\begin{equation*}
  \C{K} := \sum_{l=1}^{K+1}A_{(l)}B^{(K+2-l)}
\end{equation*}
will require that $p \ge T' + Kt - 1$, and since the spare carry locations decrease by 1 at each diagonal, diagonal $K > 0$ will have $2^{t} - K$ spare carry locations---the case $K=0$ does not require a shift and therefore diagonal 0 has no spare carry locations.
It is immediate to see that the number of spare carry locations will become negative as soon as $K > 2^{t}$.

We can continue filling up a level by accumulating partial products exactly, until we either (1) reach a diagonal $K$ such that $p < T' + Kt + 1$, or (2) hit a diagonal that contains $2^{t}$ products.
In the first case, we need to start a new level, as we have exhausted the number of adjacent diagonals that can be accumulated exactly in the current floating-point format.
If the algorithm reaches a diagonal that contains more than $2^{t} - 1$ products, then shifting the current diagonal right by $t$ bits will not be sufficient to account for the up to $t+1$ bits needed to represent the carries, and we need to check whether the additional carry bits produced can be allocated to a spare carry location of one of the preceding diagonals.

A diagonal can have at most $\min(\sA, \sB)$ products to accumulate, and this value is unlikely to exceed $2^{t}-1$ for the values of $t$ and $p$ currently of interest.
Nevertheless, we discuss how to address this situation for \cref{alg:ooy-24}, since smaller values of $t$ may become available in the future, and larger values of $\sA$ or $\sB$ might become necessary.

In order to understand whether an additional bit should be used for carries, we need to compute the maximum number of spare carry locations for a level starting with diagonal $K$.
A level starting with diagonal $K$ and ending with diagonal $M$ will have
\begin{equation}
  \label{eq:spare-carry-positions}
  \begin{aligned}
  \eta(K,M)
  &:= \sum_{l = K}^{M} (2^{t} - l) %
  = \frac{M-K+1}{2}(2^{t+1} - M - K)
  \end{aligned}
\end{equation}
spare carry positions.
As long as $\eta(K,N)$ in~\eqref{eq:spare-carry-positions} is non-negative, no additional bits are needed, and if $\eta(K,N) < 0$, then the number of extra bits needed will be $\ceil{\log_{2} -\eta(K,N)}$.
Therefore, level $j$ can accommodate
\begin{equation}
  \label{eq:def-Qj}
  Q_{j} := \left\lfloor \frac{p - T' - 1 - \ceil{\log_{2} \max\{1, -\eta(K_{j}, K_{j}+Q_{j})\})}}{t} \right\rfloor,
\end{equation}
where $K_{j}$ is the index of the first diagonal in level $j$.
Definition~\eqref{eq:def-Qj} provides the exact number of diagonals at level $j$ in the worst case, but its recursive nature makes it difficult to use in practice.
We can obtain a bound on the maximum value of $\eta(K_{j}, K_{j}+Q_{j})$ by noticing that $K_{j} \ge 0$ and $K_{j}+Q_{j} \le \max\{\sA, \sB\}$.
In fact, we can exclude the diagonal 0, which does not have any spare carry locations and only contribute one carry that is already accounted for explicitly in~\eqref{eq:def-Qj}, and look at $\eta(1, \sA)$.
If this quantity is positive, then the term reduces to $\ceil{\log_{2}0}$ and disappears; if it is negative, then evaluating the logarithm will give us a lower bound on the number of diagonals per level, which in turn will overestimate the number of inexact floating-point additions and provide an upper bound on the rounding error.

For binary64 ($p=53$), assuming an integer MMA that accepts INT8 ($t = 7$) inputs and produces INT32 ($T = 31$) outputs, we have the lower bound
\begin{equation*}
  Q_{j} \le \left\lfloor
    \frac{53-31-1}
    {7}
  \right\rfloor = 3.
\end{equation*}
Solving~\eqref{eq:spare-carry-positions} for $M$, with $K = 1$ and $t = 7$, reveals that this bound is tight as long as
$\max\{\sA,\sB\} \le 255$, which covers all practical cases of interest.

\subsection{Discussion}
\label{sec:ea-discussion}

What~\eqref{eq:bound-full} is saying is that the overall error can be substantial if either $\kA$ or $\kB$ are large.
One can counteract the prominence of these two terms by increasing $\sA$ and $\sB$, but doing so will negatively impact the performance of the algorithm, which performs $\chi(\sA,\sB)$ integer matrix multiplications.
Choosing a larger $\sA$ or $\sB$ will also increase the constant $\gamma_{\psi-1}$,
but this will only have a marginal effect on the bound, as long as $u$ is sufficiently small---especially if the accumulation is mostly done in integer arithmetic, in which case $\psi = \chi(\sA,\sB)$.

In principle, we could use $\kA$ and $\kB$ as an inexpensive way to determine $\sA$ and $\sB$: these two quantities can be computed at a negligible extra cost, as it is necessary to scan all the entries of $A$ and $B$ to compute the scale vectors $\alpha$ and $\beta$.
But what should the target value of $\kA \uA$ and $\kB \uB$ be?
In principle, we could try to choose $\sA$ and $\sB$ to minimize the product $\sA \sB$, which dictates the performance of the method, while keeping the higher-order terms approximately equal, that is,
\begin{equation}
  \label{eq:err-bound}
  \kA 2^{-\sA t} + \kB 2^{-\sB t} \approx \gamma_{\psi}.
\end{equation}
If $\sA$ and $\sB$ were reals, this would be a non-linear, constrained optimization problem in two variables with non-linear constraints.
However, we can only pick them as positive integers, and we know that for the algorithm to be efficient they cannot be too large.
Therefore, we can simply evaluate the left-hand side of~\eqref{eq:err-bound} for small values of $\sA$ and $\sB$, and take the combination that minimizes $\chi(\sA,\sB)$ among all those that deliver the correct error bound.

\section{Numerical experiments}
\label{sec:numerical-experiments}

The goal of our numerical evaluation is twofold.

On the one hand, we want to validate numerically the error analysis in~\cref{sec:error-analysis}.
We run these experiments in MATLAB using the \texttt{gemmi}\footnote{\url{https://github.com/north-numerical-computing/gemmi}} library.
This library is a flexible implementation of the Ozaki scheme, and supports all variants discussed in \cref{sec:ozaki-scheme}.
It is written in C++ but can be used in MATLAB through a complete MEX interface that exposes all functionalities.
The MATLAB source code used and the instructions on how to regenerate the data and the plots, are available.\footnote{\url{https://github.com/north-numerical-computing/integer-matrix-multiply-experiments}}

On the other hand, we want to gauge the performance of existing implementations of the Ozaki scheme on current and future hardware.
These tests are performed on the latest NVIDIA Grace-Hopper system and on a Blackwell GPU, to which we received early access.

In the experiments, we often use random matrices to illustrate our point.
We denote by $\mathcal{N}(\mu, \sigma)$ the normal distribution with mean $\mu$ and variance $\sigma^{2}$, and by $\mathcal{U}(a, b)$ the uniform distribution over the open interval $(a,
 b)$.
From both distributions we sample binary64 values.

\subsection{Behaviour on badly scaled matrices}
\label{sec:minimal-example}

\pgfplotsset{
    surfstyle/.style={
        view={0}{90},
        enlargelimits=false,
        axis on top,
        colormap/viridis,
        point meta=explicit,
        point meta min=-20,
        point meta max=0,
        xlabel={$\varphi$},
    }
}
\begin{figure}
  \centering
  \begin{tikzpicture}[trim axis group left, trim axis group right]
      \begin{groupplot}[
        group style={
          group size=2 by 1,
          horizontal sep=1.5cm,
          vertical sep = 2cm
        },
        ylabel near ticks,
        every axis plot/.append style={thick},
        ]
        \nextgroupplot[
        title={Algorithm of~\cite{ooy24}.},
        ylabel={$s$},
        surfstyle, colorbar style={
          ymode=log,
          point meta min=1e-16,
          point meta max=1e0,
          width=0.2cm,
          at={(1.1,1)}
        }]
        \addplot3[draw=none,surf,shader=interp, mesh/rows=30]
            table [meta=magnitude] {./data/fwd_err_gemmi_ooy.dat};

      \nextgroupplot[
        title={Algorithm of~\cite{uoi25}.},
        surfstyle, colorbar, colorbar style={
          ymode=log,
          point meta min=1e-16,
          point meta max=1e0,
          width=0.2cm,
          at={(1.1,1)}}
        ]
      \addplot3[draw=none,surf,shader=interp, mesh/rows=30]
          table [meta=magnitude] {./data/fwd_err_gemmi_uoi.dat};

      \end{groupplot}
  \end{tikzpicture}
\caption{Error~\eqref{eq:fwd-err-def} for the vectors in~\eqref{eq:simple-dot-product} with $\varphi$ between 0 and 100.}
\label{fig:dot-product-2}
\end{figure}
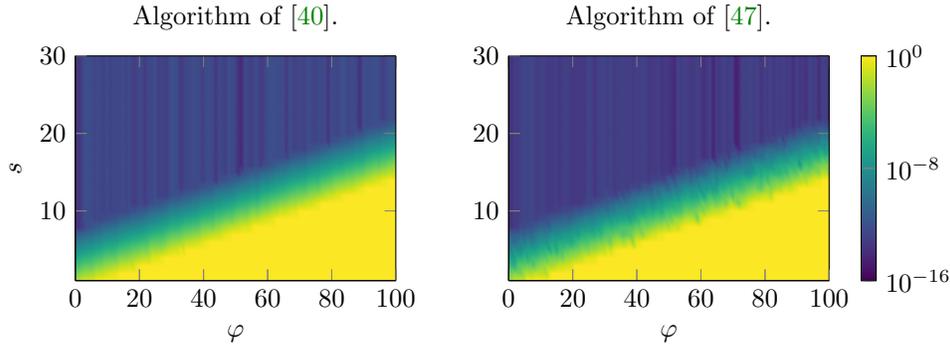

As a minimal example, we consider the computation of the inner product $a^{T}b$, where
\begin{equation}
  \label{eq:simple-dot-product}
  a = \begin{bmatrix} 2^{-\varphi}x \\ 1 \end{bmatrix},\qquad
  b = \begin{bmatrix} 2^{\varphi}y \\ 1 \end{bmatrix},\qquad
  x, y \sim \mathcal{N}(0, 1).
\end{equation}
\Cref{fig:dot-product-2} reports the relative forward error
\begin{equation}
  \label{eq:fwd-err-def}
  \frac{\abs{\wh c - c}}{\abs{c}},
\end{equation}
where $\wh c$  is the inner product $a^{T}b$ computed using a variant of the Ozaki scheme with $T=31$ and $t'=7$, and $c$ is a reference solution computed using the MATLAB Symbolic Toolbox with 32 decimal digits of accuracy.

For the vectors in~\eqref{eq:simple-dot-product}, we have that $\kA = 2^{\varphi+1} x $ and $\kB = 2^{\varphi+1}y$, and the results show that the more badly scaled the vectors are, the more slices are necessary to obtain an accurate result.
For $\varphi = 0$, about $7$ slices are sufficient to obtain binary64 accuracy, while for $\varphi = 100$ over 20 slices are needed to obtain the same accuracy.
The two variants to the Ozaki scheme considered performs similarly on this example.

Next, we extend the experimental setup in~\cite[sect.~4.2.1]{ooy24} to include matrices with a wider dynamic range.
We take the matrices $A \in \R^{10\times k}$ and $B \in \R^{k \times 10}$ with
\begin{equation}
    \label{eq:rand-mat}
    A_{ij} = a_{ij} e^{\varphi x_{ij}}, \quad
    B_{ij} = b_{ij} e^{\varphi y_{ij}}, \quad
    a_{ij}, b_{ij} \sim \mathcal{U}(-0.5, 0.5), \quad
    x_{ij}, y_{ij} \sim \mathcal{N}(0,1),
\end{equation}
where we sample from the uniform and normal distribution using the MATLAB functions \texttt{rand} and \texttt{randn}, respectively.
With the parameter $\varphi$, we can control the exponent range of the entries of $A$ and $B$ in~\eqref{eq:rand-mat}.
Previous work~\cite{ooy24} considers $\varphi \in \{0.1, 1, 2, 4\}$.

\begin{figure}[t]
  \begin{center}
    \begin{tikzpicture}[trim axis group left, trim axis group right]
      \begin{groupplot}[
        group style={
          group size=2 by 3,
          horizontal sep=1.5cm,
          vertical sep = 2cm
        },
        grid=major,
        ymode = log,
        ymax = 5e2,
        ymin=2e-19,
        xmode = log,
        ylabel near ticks,
        every axis plot/.append style={thick},
        cycle list name=test_matmul_accuracy
        ]

        \nextgroupplot[
        ylabel={Maximum error},
        title={$\varphi = 8$},
        xlabel = {$k$}
        ]
        \addplot table [x=n, y=standard-binary64] {./data/test_matmul_berr_accuracy_phi8.dat};
        \addplot table [x=n, y=split2] {./data/test_matmul_berr_accuracy_phi8.dat};
        \addplot table [x=n, y=split4] {./data/test_matmul_berr_accuracy_phi8.dat};
        \addplot table [x=n, y=split6] {./data/test_matmul_berr_accuracy_phi8.dat};
        \addplot table [x=n, y=split8] {./data/test_matmul_berr_accuracy_phi8.dat};
        \addplot table [x=n, y=split10] {./data/test_matmul_berr_accuracy_phi8.dat};

        \nextgroupplot[
        title={$\varphi = 13$},
        xlabel = {$k$},
        ]
        \addplot table [x=n, y=standard-binary64] {./data/test_matmul_berr_accuracy_phi13.dat};
        \addplot table [x=n, y=split2] {./data/test_matmul_berr_accuracy_phi13.dat};
        \addplot table [x=n, y=split4] {./data/test_matmul_berr_accuracy_phi13.dat};
        \addplot table [x=n, y=split6] {./data/test_matmul_berr_accuracy_phi13.dat};
        \addplot table [x=n, y=split8] {./data/test_matmul_berr_accuracy_phi13.dat};
        \addplot table [x=n, y=split10] {./data/test_matmul_berr_accuracy_phi13.dat};

      \end{groupplot}
    \end{tikzpicture}

    \begin{tikzpicture}[trim axis left, trim axis right]
      \begin{axis}[
        title = {},
        legend columns=3,
        scale only axis,
        width=1mm,
        height=1mm,
        hide axis,
        /tikz/every even column/.append style={column sep=0.4cm},
        legend style={at={(0,0)},anchor=center,draw=none,
          legend cell align={left},cells={line width=0.75pt}},
        legend image post style={sharp plot},
        legend cell align={left},
        cycle list name=test_matmul_accuracy
        ]
        \addplot (0,0);
        \addplot (0,0);
        \addplot (0,0);
        \addplot (0,0);
        \addplot (0,0);
        \addplot (0,0);
        \legend{binary64, $s=2$, $s=4$, $s=6$, $s=8$, $s=10$};
      \end{axis}
    \end{tikzpicture}
  \end{center}
  \caption{Error~\eqref{eq:fwd-err-def-mat} obtained by replicating the set up in \cite{ooy24}.}
  \label{fig:test_matmul_accuracy}
\end{figure}
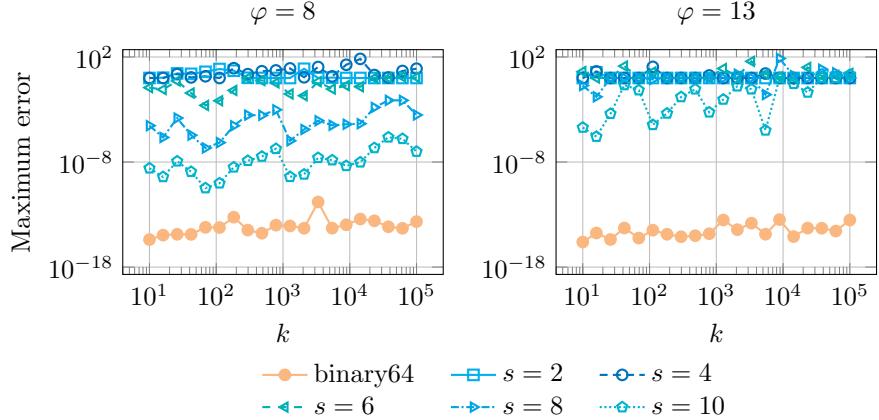

Following previous work~\cite{ooy24}, we measure the element-wide relative error
\begin{equation}
\label{eq:fwd-err-def-mat}
\max_{i,j}\frac{\abs{\widehat{c}_{ij}-c_{ij}}}{\abs{c_{ij}}},
\end{equation}
where $\wh C$ is the matrix product $AB$ computed using a variant of the Ozaki scheme, and $C$ is a reference solution computed using the MATLAB Symbolic Toolbox with 32 decimal digits of accuracy,

In~\cref{fig:test_matmul_accuracy}, we show the accuracy of the different algorithms for $\varphi=8$ (left) and $\varphi=13$ (right).
For $\varphi=13$, even the $10$-slice variant of the Ozaki scheme produces large errors.
Accuracy could be improved by using a more slices, which would reduce the truncation error in the inputs at the price of an increased runtime.
Depending on the relative performance of the integer and binary64 arithmetics available, this might make the Ozaki scheme impractical.

\pgfplotscreateplotcyclelist{test_matmul_accuracy}{
    reference style\\
    ozimmu s 2 style\\
    ozimmu s 4 style\\
    ozimmu s 6 style\\
    ozimmu s 8 style\\
    ozimmu s 10 style\\
}

\subsection{Solving $Ax=b$ with block LU factorization}
\label{sec:lu-experiment}

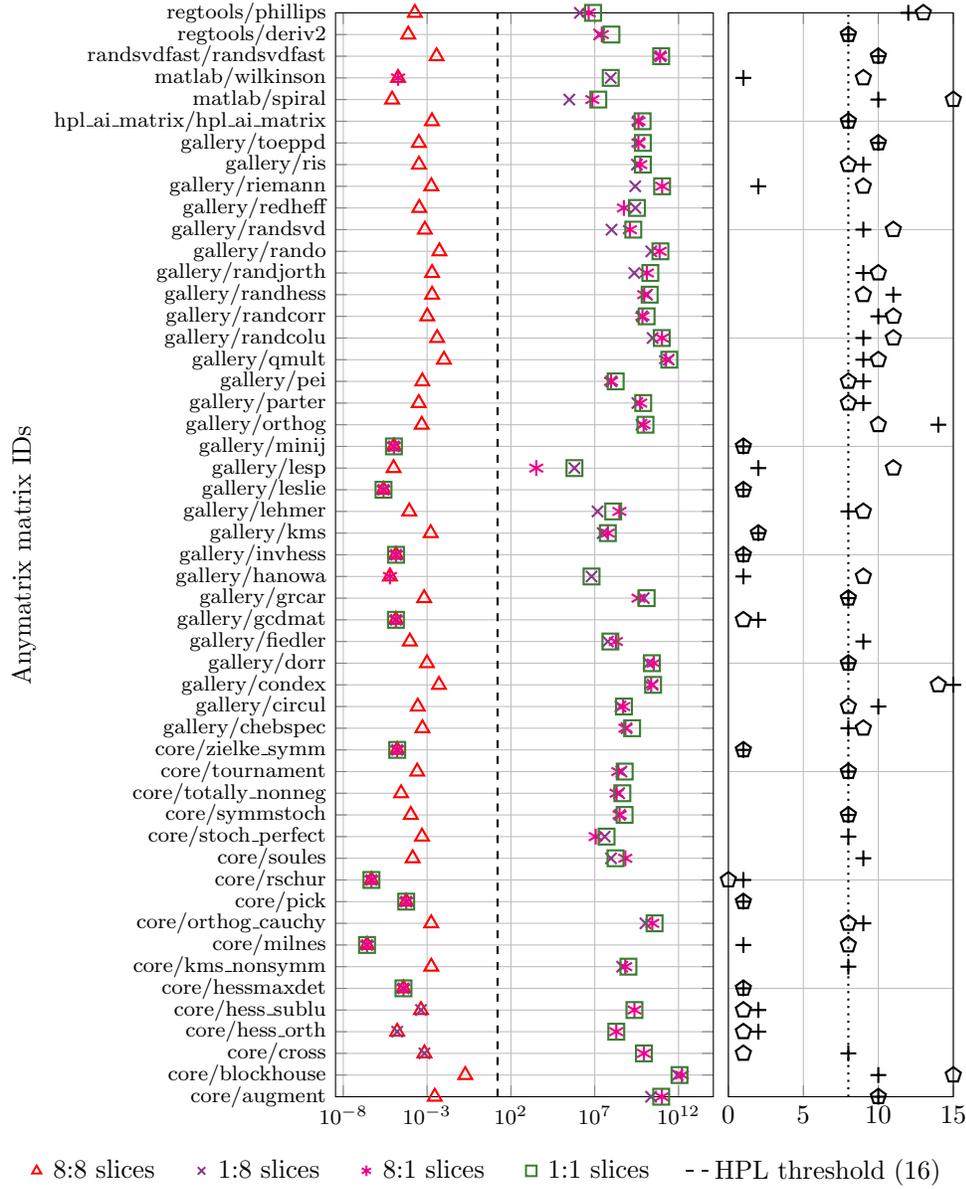
\begin{figure}[htbp]
  \begin{center}
    \begin{tikzpicture}
      \begin{groupplot}[
        group style={
          group size=2 by 1,
          horizontal sep=0.2cm,
        },
        width=2.6in,
        height=6.3in,
        grid=major,
        xlabel near ticks,
        every axis plot/.append style={thick},
        ]

        \nextgroupplot[
        title={},
        ylabel = {Anymatrix matrix IDs},
        xmode = log,
        yticklabels from table={./data/gaussian_IMMA_test.dat}{matrixID},
        ytick=data,
        yticklabel style={
          rotate=0
        },
        ticklabel style = {font=\footnotesize},
        mark options={yscale=1.5, xscale=1.5},
        ymin = 0,
        ymax = 50,
        ]

        \addplot[color=red, mark=triangle, only marks] table [y expr=\coordindex, x=88] {./data/gaussian_IMMA_test.dat};
        \addplot[color=Fuchsia, mark=x, only marks] table [y expr=\coordindex, x=18] {./data/gaussian_IMMA_test.dat};
        \addplot[color=magenta, mark=asterisk, only marks] table [y expr=\coordindex, x=81] {./data/gaussian_IMMA_test.dat};
        \addplot[color=OliveGreen, mark=square, only marks] table [y expr=\coordindex, x=11] {./data/gaussian_IMMA_test.dat};
        \addplot[dashed, black, thick] coordinates { (16,0) (16,55)};

        \nextgroupplot[
        width=1.8in,
        height=6.3in,
        ymin = 0,
        ymax = 50,
        xmin = 0,
        xmax = 15,
        ymajorticks=false,
        mark options={yscale=1.5, xscale=1.5},
        ]

        \addplot[xbar, color=black, mark=pentagon, only marks] table [y expr=\coordindex, x=splitsA] {./data/gaussian_IMMA_test.dat};
        \addplot[xbar, color=black, mark=+, only marks] table [y expr=\coordindex, x=splitsB] {./data/gaussian_IMMA_test.dat};
        \addplot[dotted, black, thick] coordinates {(8,0) (8,55)};

      \end{groupplot}
    \end{tikzpicture}\\[5pt]
    \begin{tikzpicture}[trim axis left, trim axis right]
      \begin{axis}[
        title = {},
        legend columns=5,
        scale only axis,
        width=1mm,
        height=1mm,
        hide axis,
        /tikz/every even column/.append style={column sep=0.4cm},
        legend style={at={(0,0)},anchor=center,draw=none,
          legend cell align={left}, cells={line width=0.75pt}},
        legend image code/.code={
          \draw[draw=none] (-0.15cm,0cm) -- (0.15cm,0cm);
          \draw[#1] plot coordinates {(0cm,0cm)};
        },
        legend image post style={sharp plot},
        legend cell align={left},
        ]
        \addplot [color=red, mark=triangle, only marks] (0,0);
        \addplot [color=Fuchsia, mark=x, only marks] (0,0);
        \addplot [color=magenta, mark=asterisk, only marks] (0,0);
        \addplot [color=OliveGreen, mark=square, only marks] (0,0);
        \addplot [dashed, color=Black,legend image code/.code={
             \draw[dashed] plot coordinates {(-.15cm,0cm) (.15cm,0cm)};
        }](0,0);
        \legend{8:8 slices, 1:8 slices, 8:1 slices, 1:1 slices, HPL threshold (16)};
      \end{axis}
    \end{tikzpicture}
  \end{center}
\caption{Error~\eqref{eq:hpl-measure} in solving $Ax=b$ with block LU factorization and integer-based Ozaki scheme. The right panel of the figure shows the minimum number of slices needed to split two $L_{21}$ and $U_{12}$ with no loss of information across all steps of the block LU factorization. The pentagonal and plus markers represent the minimum number of slices needed to represent exactly $L_{21}$ and $U_{12}$, respectively; the number of slices is not reported if greater than 15. The dotted line marks 8 slices.}
  \label{fig:gaussian}
\end{figure}
\begin{figure}
  \centering
  \begin{tikzpicture}
    \begin{axis}[
      view={0}{90},
      colorbar,
      colorbar style={
        ytick={-2,3,8,13},
        yticklabels={$10^{-2}$,$10^3$,$10^{8}$,$10^{13}$},
        ylabel={Quantity in~\eqref{eq:hpl-measure}},
        width=0.2cm,
      },
      xlabel={$s_{A}$},
      ylabel={$s_{B}$},
      colormap/viridis,
      point meta=log10(\thisrow{rel_error}),
      point meta min=-2,
      point meta max=13,
      xtick distance=1,
      ytick distance=1,
      height=15cm,
      axis equal image,
      xmin=0.5, xmax=8.5,
      ymin=0.5, ymax=8.5,
      tick style={draw=none},
      axis line style={draw=none},
      clip=false
      ]

      \addplot[
      matrix plot,
      mesh/cols=8,
      ]
      table[
      x=splitA,
      y=splitB,
      meta=rel_error
      ]{./data/gaussian_IMMA_test_exhaustive_splits.dat};

      \draw[BurntOrange!40, ultra thick, densely dotted] (6.5,6.5) rectangle (8.5,8.5);
    \end{axis}
  \end{tikzpicture}

  \caption{Error~\eqref{eq:hpl-measure} in solving $Ax=b$ with block LU factorization, where $A$ is the Anymatrix matrix \texttt{core/blockhouse} with $n=500$.
    At each step of the factorization, the algorithm multiplies the matrices
    $L_{21}$ and $U_{12}$ using the integer-based Ozaki scheme with $1$ to $8$
    slices.
    The error is below the HPL threshold (16) for configurations within the
    dotted square (seven to eight slices for both matrices).}
  \label{fig:gaussian-one-mat}

\end{figure}

We now investigate the performance of the Ozaki scheme on higher-level matrix algorithms.
The computers in the TOP500 list are ranked based on their performance on the HPL benchmark, which solves a linear system of equations using LU factorization~\cite{dlp03},~\cite{pwdc18}.
We study the accuracy of an implementation of the block LU factorization~\cite[sect.~3.6.1]{gova13} where the update in the Schur complement is computed using the Ozaki scheme.

In our experiments, we consider the linear system $Ax=b$ of order $n=500$, where $A \in \R^{n\times n}$ is one of the nonsingular test matrices from the anymatrix matrix collection~\cite{himi21},~\cite{himi21-UG} and $b \in \R^{n}$ has entries sampled from $\mathcal U(0,1)$.
$A$, $x$, and $b$ are represented in binary64.
We chose a subset of anymatrix matrices for which the solution was not exact and for which MATLAB did not produce a warning about the test matrix being close to singular.
The solution $x$ is computed using block LU factorization with partial pivoting and block size $b = 10$.

Let $A^{(i-1)}$ be the trailing submatrix at step $i$ of the block LU factorization.
Consider the partitioning
\begin{equation*}
A^{(i-1)} =:
\begin{bmatrix}
A_{11} & A_{12} \\
A_{21} & A_{22}
\end{bmatrix} \in \R^{(n - b(i-1)) \times (n - b(i-1))},
\end{equation*}
where $A_{11} \in \R^{b \times b}$ is the current panel, while $A_{12} \in \R^{b \times (n - bi)}$, $A_{21} \in \R^{(n - bi)\times b}$, and $A_{22} \in \R^{(n - bi) \times (n - bi)}$.
The panel factorization produces
\begin{equation*}
    A_{11} =: L_{11}U_{11}, \qquad
    U_{12} = L_{11}^{-1} A_{12}, \quad
    L_{21} = A_{21} U_{11}^{-1},
\end{equation*}
and the trailing sumbatrix $A_{22}$ is updated via the Schur complement
\begin{equation}
    \label{eq:schur-complement}
  A^{(i)} := A_{22} - L_{21} U_{12}.
\end{equation}
In our experiments, we compute the outer product in~\eqref{eq:schur-complement} using the Ozaki scheme with $T=31$ and $t'=7$, and with four different splitting configurations: eight slices per matrix, one slice per matrix, and one slice for one matrix and eight for the other.

For each matrix in our test set, we report the backward error of the computed solution~$\wh x$ in terms of the measure used in the HPL benchmark
\begin{equation}
  \label{eq:hpl-measure}
  \frac{\|A \wh x-b\|_\infty}{2u (\|A\|_{\infty}\|\wh x\|_{\infty}+\|b\|_\infty) n}.
\end{equation}
An HPL run is declared failed if the quantity in~\eqref{eq:hpl-measure} is above the threshold value 16.

To understand how many slices are needed throughout the algorithm, we also compute, at each step $i$, the minimum number of slices needed to represent $L_{21}$ and $U_{12}$ exactly---in the language of \cref{sec:ea-truncation-error}, we compute $\sAs$ for $L_{21}$ and $\sBs$ for $U_{12}$.
In \cref{fig:gaussian}, we report the maximum across all steps of the block LU factorization.

The minimum number of slices that preserves all the bits in the original floating-point matrix is computed as follows.
First, we calculate the \emph{bit spread}, the number of bits needed to represent the significand of each entry of $L_{21}$ without truncation error.
This will be the number of bits between the most significant and the least significant bit of the significand set to 1, inclusive.
For example, if the significand of the entry only contains zeros, then bit spread will be zero, and if it only has a single bit set to one, then the bit spread will be one.
In general, the bit spread will be an integer between zero and 53 in binary64.
For each row of $L_{21}$, we also calculate the difference between the largest and smallest exponents, and add it to the maximum element bit spread for that row.
We take the maximum value thus obtained across all rows, and we divide it by 7 since $t'=7$---this is $\sAs$, the minimum number of slices needed to represent $L_{21}$ exactly.
We repeat the same process column-wise on $U_{12}$ to obtain $\sBs$.

\cref{fig:gaussian} shows that the number of slices needed to meet the HPL threshold for the measure~\eqref{eq:hpl-measure} is matrix dependent.
For $\sA = \sB = 8$, all matrices in our test set meet the threshold.
Eight of the 51 matrices (for example, \texttt{gallery/minij}) meet the threshold regardless of the number of slices used.
A subset of matrices displays a dependency on specific values of $\sA$ or $\sB$.
For \texttt{matlab/wilkinson} and \texttt{gallery/hanowa}, $\sAs=9$ and $\sBs=1$, and these two matrices satisfy the threshold with $\sA = 8$ and $\sB = 1$, but not vice versa.
Conversely, for \texttt{core/hess\_orth} and \texttt{gallery/hess\_sublu} we have $\sAs=1$ and $\sBs=2$, and for \texttt{core/cross} we have $\sAs=1$ and $\sBs=8$. These matrices satisfy the threshold with $\sA = 1$ and $\sB = 8$.

\Cref{fig:gaussian-one-mat} demonstrates how the error \eqref{eq:hpl-measure} varies by changing the number of slices used by the Ozaki scheme between 1 and 8.
The error is below the HPL threshold of 16 when at least 7 slices are used for both matrices.

These results support the conclusions in \cref{sec:error-analysis} that, for some inputs, the number of slices used in the Ozaki scheme can be reduced without any loss of accuracy.

\subsection{Benchmarking on GPUs}
\label{sec:benchmarking-gpus}
\newcommand{\dgemm}{\texttt{DGEMM}}
\newcommand{\cublas}{\texttt{cublas}}
\newcommand{\ozimmu}[1][]{\texttt{ozIMMU#1}}
\newcommand{\cuimma}[1][]{\texttt{cuIMMA#1}}

We now examine the performance of the algorithms by benchmarking them on NVIDIA GPUs.
A general matrix multiplication (GEMM) operation has the form $D = \alpha A B + \beta C$, where $A \in \R^{m \times k}$, $B \in \R^{k \times n}$, and $C,D \in \R^{m \times n}$. The standard implementation of GEMM requires $2mnk$ floating-point operations.
In our experiments, we focus on the BLAS routine \dgemm{}, which computes a GEMM where all matrices and scalars are binary64 values.
We compare two implementations.

\begin{itemize}
  \item \cublas{} is the \dgemm{} implementation in the cuBLAS library, which uses binary64 arithmetic throughout.
  \item \cuimma{} is developed by NVIDIA as part of the cuBLAS library---we gained early access to this new development in the form of a prototype that is not yet publicly available. The engineering cuBLAS build refers to this new feature as ``FP64 Emulation through the IMMA instructions''. We do not know what version of the algorithm is implemented by this engineering build, but we note NVIDIA has recently published some details about their implementation of this algorithm~\cite{sabb26}.
\end{itemize}

We do not consider the implementations by Ootomo\footnote{\url{https://github.com/enp1s0/ozIMMU}} and Uchino\footnote{\url{https://github.com/RIKEN-RCCS/accelerator\_for\_ozIMMU}} because, in the experiments on badly scaled matrices at the end of this section, they exhibit some unexpectedly large errors that suggest these codes might be unstable in some cases.

\cuimma{} requires a fairly large workspace attached to the cuBLAS handle, which can be set with the \texttt{cublasSetWorkspace} function. The code has been integrated into the MAGMA library~\cite{abcg24,addh09} to leverage MAGMA's testing capabilities and to test them within higher-level LAPACK algorithms on the GPU.

The tests are performed on two NVIDIA GPUs: a Grace-Hopper system (GH200), which features a 72-core Grace CPU and an H100 GPU (released in 2022) with 96GB of HBM3e memory, and a production-level Blackwell B200 GPU (released in 2025), equipped with approximately 178GB of HBM3e memory. Two separate instances of MAGMA were compiled to use CUDA 12.2 on the GH200 system and a prototype of CUDA 12.8 with the B200~GPU. MAGMA also requires a CPU LAPACK library, which is used as reference implementation to assess the accuracy of the algorithms. We use the NVIDIA Performance Libraries (\texttt{NVPL}) on the GH200 system and the Intel \texttt{MKL} Library with the B200~GPU.
Because the behavior of the numerical accuracy is very similar for the two GPUs, when both results are available we only show the accuracy on the B200~GPU.

\paragraph{Benchmarks using random matrices}
We begin by running the experiments on random matrices generated with \texttt{DLARNV}, which have entries sampled from $\mathcal U(0,1)$. The purpose of these benchmarks is to observe the performance and accuracy of square matrix multiplication, which is usually a good benchmark for measuring the asymptotic performance on GPUs.

For the accuracy, we use the forward error
\begin{equation}
\frac{\left\| \wh D-D \right\|_{F}}{|\alpha|\sqrt{k+2}\left\| A \right\|_{F}\left\| B \right\|_{F} + 2|\beta| \left\| C \right\|_{F}},
\label{eq:err_higham}
\end{equation}
where $\wh D$ is the result computed on the GPU and $D$ is a reference result computed by the BLAS implementation running on the CPU.
This measure, based on~\cite[sect.~3.5]{high02}, is used in the MAGMA testing suite for matrix multiply, where a test is considered ``passed'' if the quantity~\eqref{eq:err_higham} is below unit roundoff.

\pgfplotscreateplotcyclelist{dgemm_perf_sq}{
    reference style\\
    cuimma s 3 style\\
    cuimma s 6 style\\
    cuimma s 7 style\\
}

\begin{figure}[t]
  \begin{center}
    \begin{tikzpicture}[trim axis group left, trim axis group right]
      \begin{groupplot}[
        group style={
          group size=2 by 1,
          horizontal sep=1.5cm,
          vertical sep = 2cm
        },
        grid=major,
        ylabel near ticks,
        every axis plot/.append style={thick},
        ymax=450000,
        ymin=-50000,
        xmax=22000,
        xmin=-2000,
        cycle list name=dgemm_perf_sq
        ]

        \nextgroupplot[
        ylabel={Gflop/s},
        title={GH200},
        xlabel = {$m=k=n$}
        ]
        \addplot table [x=m, y=cublas] {./data/GPU_raw_data/GH200/dgemm_imma/dgemm_sq_ozimmu_summary_gh200.txt};
        \addplot table [x=m, y=oz03-perf] {./data/GPU_raw_data/GH200/dgemm_imma/dgemm_sq_cuimma_summary_gh200.txt};
        \addplot table [x=m, y=oz06-perf] {./data/GPU_raw_data/GH200/dgemm_imma/dgemm_sq_cuimma_summary_gh200.txt};
        \addplot table [x=m, y=oz07-perf] {./data/GPU_raw_data/GH200/dgemm_imma/dgemm_sq_cuimma_summary_gh200.txt};

        \nextgroupplot[
        title={B200},
        xlabel = {$m=k=n$}
        ]
        \addplot table [x=M, y=cuBLAS_gflops] {./data/GPU_raw_data/B200/dgemm_imma/dgemm_sq_cuimma_s3.txt};

        \addplot table [x=M, y=MAGMA_gflops] {./data/GPU_raw_data/B200/dgemm_imma/dgemm_sq_cuimma_s3.txt};
        \addplot table [x=M, y=MAGMA_gflops] {./data/GPU_raw_data/B200/dgemm_imma/dgemm_sq_cuimma_s6.txt};
        \addplot table [x=M, y=MAGMA_gflops] {./data/GPU_raw_data/B200/dgemm_imma/dgemm_sq_cuimma_s7.txt};

      \end{groupplot}
    \end{tikzpicture}

    \begin{tikzpicture}[trim axis left, trim axis right]
      \begin{axis}[
        title = {},
        legend columns=4,
        scale only axis,
        width=1mm,
        height=1mm,
        hide axis,
        /tikz/every even column/.append style={column sep=0.4cm},
        legend style={at={(0,0)},anchor=center,draw=none,
          legend cell align={left},cells={line width=0.75pt}},
        legend image post style={sharp plot},
        legend cell align={left},
        cycle list name=dgemm_perf_sq
        ]
        \addplot  (0,0);
        \addplot  (0,0);
        \addplot  (0,0);
        \addplot  (0,0);
        \legend{\cublas{},
                \cuimma{} $s=3$, \cuimma{} $s=6$, \cuimma{} $s=7$};
      \end{axis}
    \end{tikzpicture}
  \end{center}
  \caption{Performance of \cublas{} and \cuimma{} for random square matrices.}
  \label{fig:dgemm_perf_sq}
\end{figure}
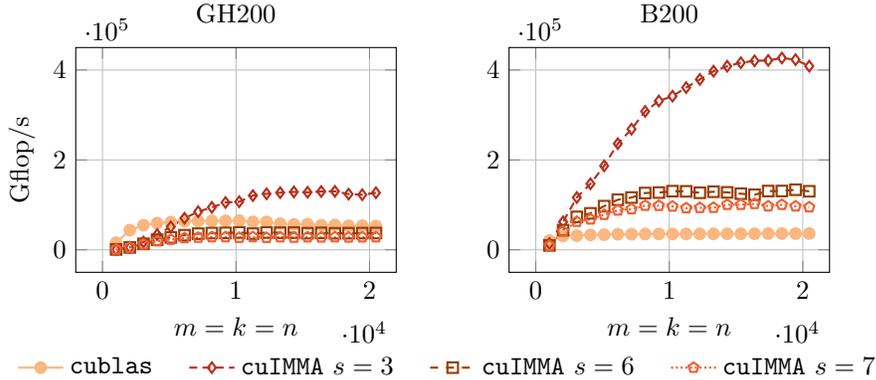

\Cref{fig:dgemm_perf_sq} compares the performance of \cublas{} with that of \cuimma{} with 3, 6, and 7 slices. A relatively small number of slices ($s = 3$) shows the best possible performance, and \cuimma{} significantly outperforms \cublas{}: the asymptotic speedup is around $2.4\times$ on the GH200 system and $7.6\times$ on the B200 GPU.
However, these huge performance gains come at a significant loss of accuracy, as shown in \cref{fig:dgemm_error_sq}. The forward error of is of order $10^{-10}$, compared with an error of order $10^{-19}$ for~\cublas{}.

\begin{figure}[t]
  \begin{center}
    \begin{tikzpicture}[trim axis group left, trim axis group right]
      \begin{groupplot}[
        group style={
          group size=2 by 1,
          horizontal sep=1.5cm,
          vertical sep = 2cm
        },
        grid=major,
        ylabel near ticks,
        every axis plot/.append style={thick},
        ymax=5e-5,
        ymin=2e-22,
        xmin=-2000,
        xmax=22000,
        ymode=log,
        cycle list name=dgemm_perf_sq
        ]

        \nextgroupplot[
        ylabel={Forward error},
        title={GH200},
        xlabel = {$m=k=n$}
        ]
        \addplot table [x=m, y=cu-err] {./data/GPU_raw_data/GH200/dgemm_imma/dgemm_sq_ozimmu_summary_gh200.txt};
        \addplot table [x=m, y=oz03-err] {./data/GPU_raw_data/GH200/dgemm_imma/dgemm_sq_cuimma_summary_gh200.txt};
        \addplot table [x=m, y=oz06-err] {./data/GPU_raw_data/GH200/dgemm_imma/dgemm_sq_cuimma_summary_gh200.txt};
        \addplot table [x=m, y=oz07-err] {./data/GPU_raw_data/GH200/dgemm_imma/dgemm_sq_cuimma_summary_gh200.txt};

        \nextgroupplot[
        title={B200},
        xlabel = {$m=k=n$}
        ]
        \addplot table [x=M, y=cuBLAS_error] {./data/GPU_raw_data/B200/dgemm_imma/dgemm_sq_cuimma_s3.txt};
        \addplot table [x=M, y=MAGMA_error] {./data/GPU_raw_data/B200/dgemm_imma/dgemm_sq_cuimma_s3.txt};
        \addplot table [x=M, y=MAGMA_error] {./data/GPU_raw_data/B200/dgemm_imma/dgemm_sq_cuimma_s6.txt};
        \addplot table [x=M, y=MAGMA_error] {./data/GPU_raw_data/B200/dgemm_imma/dgemm_sq_cuimma_s7.txt};
      \end{groupplot}
    \end{tikzpicture}

    \begin{tikzpicture}[trim axis left, trim axis right]
      \begin{axis}[
        title = {},
        legend columns=4,
        scale only axis,
        width=1mm,
        height=1mm,
        hide axis,
        /tikz/every even column/.append style={column sep=0.4cm},
        legend style={at={(0,0)},anchor=center,draw=none,
          legend cell align={left},cells={line width=0.75pt}},
        legend image post style={sharp plot},
        legend cell align={left},
        cycle list name=dgemm_perf_sq
        ]
        \addplot  (0,0);
        \addplot  (0,0);
        \addplot  (0,0);
        \addplot  (0,0);
        \legend{\cublas{},
                \cuimma{} $s=3$, \cuimma{} $s=6$, \cuimma{} $s=7$};
      \end{axis}
    \end{tikzpicture}
  \end{center}
  \caption{Error~\eqref{eq:err_higham} of \cublas{} and \cuimma{} for random square matrices.}
  \label{fig:dgemm_error_sq}
\end{figure}

Increasing the number of slices to six or seven significantly reduces the performance of the emulated \dgemm{} but improves the forward error to acceptable accuracy levels, and with seven slices \cuimma{}'s accuracy indistinguishable from that of~\cublas{}.

\cref{fig:dgemm_error_sq} also shows the promising potential of the Ozaki scheme on hardware with large performance ratios between the peak performances of binary64 and INT8. On the GH200 system, the INT8-to-binary64 performance ratio for the tensor cores is about $30$, and a 7-slice emulated \dgemm{} is slower than the binary64 implementation. In order to observe a performance gain with an acceptable accuracy, there must be a very large ratio between the theoretical peak performances of INT8 and binary64---this is the case of the B200 GPU, whose INT8-to-binary64 performance ratio is $112.5$ for the tensor cores. There, the emulated \dgemm{} implementation can outperform the floating-point implementation while maintaining a similar accuracy.

Finally, it is worth noting that \cublas{} is slower on the B200 compared with the GH200. This performance drop stems from a strategic shift in hardware design priorities, which was driven by the increasing demand for computational power at relatively low precisions, particularly for training large-scale AI models. To meet this demand, the vendor has chosen to allocate less silicon to binary64 arithmetic, focusing on enhancing support for low-precision formats.
This shift reflects a broader industry trend where hardware advancements are no longer aimed at providing uniform performance improvements across all compute precisions. Historically, new generations of hardware architectures typically brought performance gains at all precision levels, albeit at varying degrees. With the rise of deep learning and AI workloads, however, the emphasis has shifted toward optimizing for lower precisions, where performance and efficiency gains are maximized for the most needed use cases.

\paragraph{Impact on the accuracy of higher-level algorithms}
The emulated \dgemm{} routines were integrated into MAGMA's LAPACK algorithms in order to evaluate their performance in a higher-level matrix algorithm. In particular, we considered the numerical behavior of the QR factorization (\texttt{DGEQRF}), the symmetric eigensolver (\texttt{DSYEVD}), and the singular value decomposition (\texttt{DGESVD}). For all these algorithms, an acceptable accuracy was achieved when the emulated \dgemm{} used 8 rather than 7 slices. As an example, we show the accuracy of modified version of \texttt{DSYEVD} that internally uses \cuimma{} to perform matrix--matrix multiplications.

\pgfplotscreateplotcyclelist{dsyevd_error}{
    reference style\\
    cuimma s 3 style\\
    cuimma s 7 style\\
    cuimma s 8 style\\
}

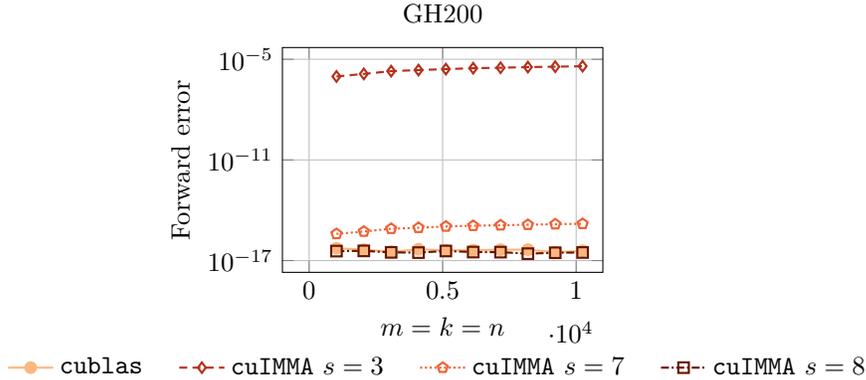
\begin{figure}[t]
  \begin{center}
    \begin{tikzpicture}[trim axis group left, trim axis group right]
      \begin{groupplot}[
        group style={
          group size=1 by 1,
          horizontal sep=1.5cm,
          vertical sep = 2cm
        },
        grid=major,
        ylabel near ticks,
        every axis plot/.append style={thick},
        ymax=5e-5,
        ymin=2e-18,
        ymode=log,
        xmax=11000,
        xmin=-1000,
        cycle list name=dsyevd_error
        ]

        \nextgroupplot[
        ylabel={Forward error},
        title={GH200},
        xlabel = {$m=k=n$}
        ]
        \addplot table [x=N, y=|I-U^H_U|] {./data/GPU_raw_data/GH200/dsyevd_imma/dsyevd_ref_poev_cluster1_1e10.txt};
        \addplot table [x=N, y=|I-U^H_U|] {./data/GPU_raw_data/GH200/dsyevd_imma/dsyevd_poev_cluster1_1e10_cuimma_s3_gh200.txt};
        \addplot table [x=N, y=|I-U^H_U|] {./data/GPU_raw_data/GH200/dsyevd_imma/dsyevd_poev_cluster1_1e10_cuimma_s7_gh200.txt};
        \addplot table [x=N, y=|I-U^H_U|] {./data/GPU_raw_data/GH200/dsyevd_imma/dsyevd_poev_cluster1_1e10_cuimma_s8_gh200.txt};
      \end{groupplot}
    \end{tikzpicture}

    \begin{tikzpicture}[trim axis left, trim axis right]
      \begin{axis}[
        title = {},
        legend columns=4,
        scale only axis,
        width=1mm,
        height=1mm,
        hide axis,
        /tikz/every even column/.append style={column sep=0.4cm},
        legend style={at={(0,0)},anchor=center,draw=none,
          legend cell align={left},cells={line width=0.75pt}},
        legend image post style={sharp plot},
        legend cell align={left},
        cycle list name=dsyevd_error
        ]
        \addplot  (0,0);
        \addplot  (0,0);
        \addplot  (0,0);
        \addplot  (0,0);
        \legend{\cublas{},
                \cuimma{} $s=3$, \cuimma{} $s=7$, \cuimma{} $s=8$};
      \end{axis}
    \end{tikzpicture}
  \end{center}
  \caption{Orthogonality~\eqref{eq:orth-eigvecs} of the eigenvectors computed by
the MAGMA symmetric eigensolver.}
  \label{fig:dsyevd_error}
\end{figure}

The benchmarks use the symmetric matrix $A = U \varLambda U^{T}$ where $U$ is an orthogonal matrix generated by applying random Householder reflectors of decreasing size, and $\lambda_i = \{1, 1, ..., 1, 10^{-10}\}$. We use $s\in\{3, 6, 7, 8\}$ and evaluate the loss of orthogonality
\begin{equation}
  \label{eq:orth-eigvecs}
  \frac{\| I - \wh U \wh U^{T}\|_{1}}{n},
\end{equation}
where $\wh U$ is the computed matrix of eigenvectors.
The results for the GH200 system are summarized in \cref{fig:dsyevd_error}.

In this experiment, \cuimma{} with $s = 8$ achieves the same level of accuracy as \cublas{}.
The results in \cref{fig:dsyevd_error} suggests that the choice of $s$, which is essential for the stability of higher-level algorithms that rely on matrix multiplication, cannot be deduced from benchmarking \dgemm{} in isolation. This conclusion could pose a challenge to the wide acceptance of emulated \dgemm{} in scientific computing, as the number of slices would have to be decided on an algorithm-by-algorithm basis.

\paragraph{Accuracy of emulated DGEMM on badly scaled matrices}
Multiplying badly-scaled matrices further exposes the weaknesses of emulated \dgemm{} algorithms.
In this test, we generate two random square matrices $A = \bar{A}D$ and $B = D^{-1}\bar{B}$, where the entries of $\bar{A}$ and $\bar{B}$ are sampled from $\mathcal U(1,2)$ and $D$ is a diagonal matrix.
This matrix is such that $d_{11} = \sqrt{\kappa_D}^{-1}$ and $d_{nn} = \sqrt{\kappa_D}$, where $\kappa_{D}$ is a parameter that controls the dynamic range of $A$ and $B$ and the ratio $d_{ii}/d_{i+1,i+1}$ is constant for $i = 1, 2, \ldots, n-1$.
Next, for $i = 1, 2, \ldots, n$ we perform a circular rotation of the $i$th row of $A$ and the $i$th column of $B$ by $i$ places. This rotation is optional but makes the benchmark more challenging.
In order to only evaluate the matrix multiplication, we set $\alpha=1$ and $\beta = 0$, which effectively simplifies the GEMM to $D = AB$.
The tests are conducted for $\kappa_D\in\{10^{10}, 10^{20}, 10^{30}, 10^{40}\}$ and $s\in\{8, 12, 16, 18\}$

For this particular test, we notice that the original formula in \eqref{eq:err_higham} could be misleading, because the product $\left\| A \right\|_{F}\left\| B \right\|_{F}$ may be very large because of the scaling described above. This is why we prefer the maximum element-wise relative error~\eqref{eq:fwd-err-def-mat}, where $C$ is the reference solution computed using BLAS on the CPU.

\begin{figure}[t]
  \begin{center}
    \begin{tikzpicture}[trim axis group left, trim axis group right]
      \begin{groupplot}[
        group style={
          group size=2 by 3,
          horizontal sep=1.5cm,
          vertical sep = 2cm
        },
        grid=major,
        ymode = log,
        ymax = 5e1,
        ymin=2e-18,
        ylabel near ticks,
        xmin=-2000,
        xmax=22000,
        every axis plot/.append style={thick},
        cycle list name=badly_scaled_dgemm_cyclist
        ]

        \nextgroupplot[
        ylabel={Maximum error},
        title={$\kappa_{D} = 10^{10}$},
        xlabel = {$m=k=n$}
        ]
        \addplot table [x=M, y=cuBLAS_error] {./data/GPU_raw_data/B200/dgemm_imma_badly_scaled/dgemm_sq_badly_scaled_Dcond_1e10_cuimma_s8_b200.txt};
        \addplot table [x=M, y=MAGMA_error] {./data/GPU_raw_data/B200/dgemm_imma_badly_scaled/dgemm_sq_badly_scaled_Dcond_1e10_cuimma_s8_b200.txt};
        \addplot table [x=M, y=MAGMA_error] {./data/GPU_raw_data/B200/dgemm_imma_badly_scaled/dgemm_sq_badly_scaled_Dcond_1e10_cuimma_s12_b200.txt};
        \addplot table [x=M, y=MAGMA_error] {./data/GPU_raw_data/B200/dgemm_imma_badly_scaled/dgemm_sq_badly_scaled_Dcond_1e10_cuimma_s16_b200.txt};
        \addplot table [x=M, y=MAGMA_error] {./data/GPU_raw_data/B200/dgemm_imma_badly_scaled/dgemm_sq_badly_scaled_Dcond_1e10_cuimma_s18_b200.txt};

        \nextgroupplot[
        title={$\kappa_{D} = 10^{20}$},
        xlabel = {$m=k=n$}
        ]
        \addplot table [x=M, y=cuBLAS_error] {./data/GPU_raw_data/B200/dgemm_imma_badly_scaled/dgemm_sq_badly_scaled_Dcond_1e20_cuimma_s8_b200.txt};
        \addplot table [x=M, y=MAGMA_error] {./data/GPU_raw_data/B200/dgemm_imma_badly_scaled/dgemm_sq_badly_scaled_Dcond_1e20_cuimma_s8_b200.txt};
        \addplot table [x=M, y=MAGMA_error] {./data/GPU_raw_data/B200/dgemm_imma_badly_scaled/dgemm_sq_badly_scaled_Dcond_1e20_cuimma_s12_b200.txt};
        \addplot table [x=M, y=MAGMA_error] {./data/GPU_raw_data/B200/dgemm_imma_badly_scaled/dgemm_sq_badly_scaled_Dcond_1e20_cuimma_s16_b200.txt};
        \addplot table [x=M, y=MAGMA_error] {./data/GPU_raw_data/B200/dgemm_imma_badly_scaled/dgemm_sq_badly_scaled_Dcond_1e20_cuimma_s18_b200.txt};

        \nextgroupplot[
        ylabel={Maximum error},
        title={$\kappa_{D} = 10^{30}$},
        xlabel = {$m=k=n$}
        ]
        \addplot table [x=M, y=cuBLAS_error] {./data/GPU_raw_data/B200/dgemm_imma_badly_scaled/dgemm_sq_badly_scaled_Dcond_1e30_cuimma_s8_b200.txt};
        \addplot table [x=M, y=MAGMA_error] {./data/GPU_raw_data/B200/dgemm_imma_badly_scaled/dgemm_sq_badly_scaled_Dcond_1e30_cuimma_s8_b200.txt};
        \addplot table [x=M, y=MAGMA_error] {./data/GPU_raw_data/B200/dgemm_imma_badly_scaled/dgemm_sq_badly_scaled_Dcond_1e30_cuimma_s12_b200.txt};
        \addplot table [x=M, y=MAGMA_error] {./data/GPU_raw_data/B200/dgemm_imma_badly_scaled/dgemm_sq_badly_scaled_Dcond_1e30_cuimma_s16_b200.txt};
        \addplot table [x=M, y=MAGMA_error] {./data/GPU_raw_data/B200/dgemm_imma_badly_scaled/dgemm_sq_badly_scaled_Dcond_1e30_cuimma_s18_b200.txt};

        \nextgroupplot[
        title={$\kappa_{D} = 10^{40}$},
        xlabel = {$m=k=n$}
        ]
        \addplot table [x=M, y=cuBLAS_error] {./data/GPU_raw_data/B200/dgemm_imma_badly_scaled/dgemm_sq_badly_scaled_Dcond_1e40_cuimma_s8_b200.txt};
        \addplot table [x=M, y=MAGMA_error] {./data/GPU_raw_data/B200/dgemm_imma_badly_scaled/dgemm_sq_badly_scaled_Dcond_1e40_cuimma_s8_b200.txt};
        \addplot table [x=M, y=MAGMA_error] {./data/GPU_raw_data/B200/dgemm_imma_badly_scaled/dgemm_sq_badly_scaled_Dcond_1e40_cuimma_s12_b200.txt};
        \addplot table [x=M, y=MAGMA_error] {./data/GPU_raw_data/B200/dgemm_imma_badly_scaled/dgemm_sq_badly_scaled_Dcond_1e40_cuimma_s16_b200.txt};
        \addplot table [x=M, y=MAGMA_error] {./data/GPU_raw_data/B200/dgemm_imma_badly_scaled/dgemm_sq_badly_scaled_Dcond_1e40_cuimma_s18_b200.txt};
      \end{groupplot}
    \end{tikzpicture}

    \begin{tikzpicture}[trim axis left, trim axis right]
      \begin{axis}[
        title = {},
        legend columns=3,
        scale only axis,
        width=1mm,
        height=1mm,
        hide axis,
        /tikz/every even column/.append style={column sep=0.4cm},
        legend style={at={(0,0)},anchor=center,draw=none,
          legend cell align={left},cells={line width=0.75pt}},
        legend image post style={sharp plot},
        legend cell align={left},
        cycle list name=badly_scaled_dgemm_cyclist
        ]
        \addplot (0,0);
        \addplot (0,0);
        \addplot (0,0);
        \addlegendimage{empty legend}
        \addplot (0,0);
        \addplot (0,0);
        \legend{\cublas{}, \cuimma{} $s=8$, \cuimma{} $s=12$, \phantom{ }, \cuimma{} $s=16$, \cuimma{} $s=18$};
      \end{axis}
    \end{tikzpicture}
  \end{center}
\caption{Error~\eqref{eq:fwd-err-def-mat} of \cublas{} and \cuimma{} on B200 GPU for badly scaled matrices.}
\label{fig:dgemm_cuimma_error_scaled}
\end{figure}
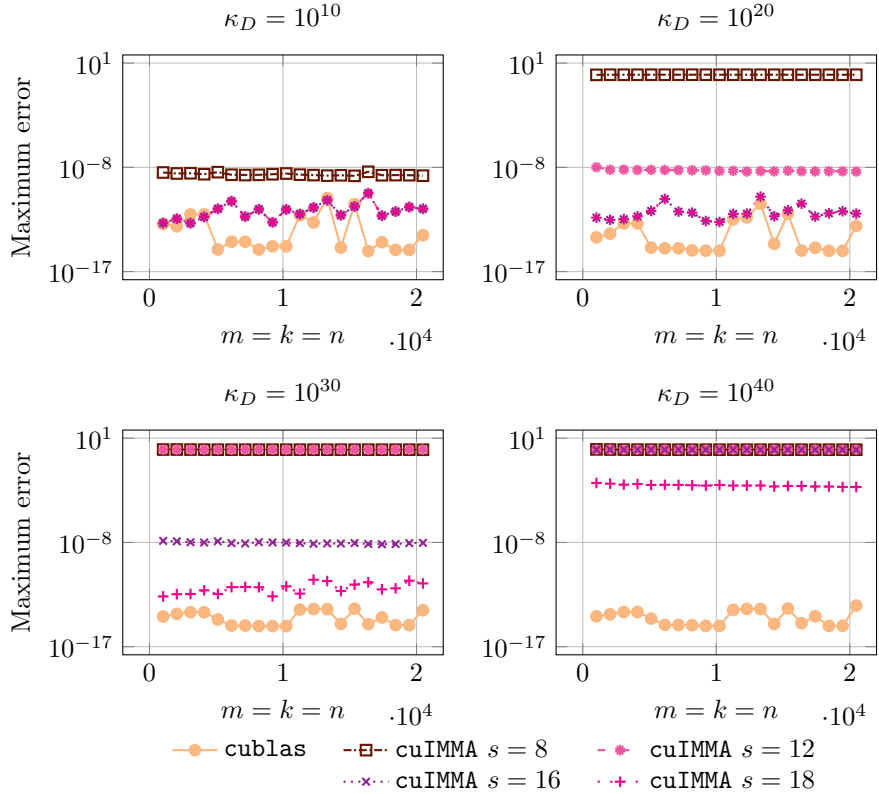

\Cref{fig:dgemm_cuimma_error_scaled} shows the results for \cuimma{}.
As predicted by the analysis, setting $s=8$, which is a safe configuration for random matrices, is not sufficient for badly scaled inputs, even for the smallest value of $\kappa_D$ considered. For lower values of $\kappa_D$, it is still possible to achieve binary64 accuracy by increasing the number of slices, but this would severely impact the performance.
For $\kappa_D\in\{10^{30}, 10^{40}\}$, the implementation becomes too inaccurate even when for $s = 18$.

This illustrates the need for a method to estimate how many slices the Ozaki scheme will require in order to achieve a prescribed level of accuracy.

\pgfplotscreateplotcyclelist{badly_scaled_dgemm_cyclist}{
    reference style\\
    cuimma s 8 style\\
    cuimma s 12 style\\
    cuimma s 16 style\\
    cuimma s 18 style\\
}

\section{Conclusion}

The Ozaki scheme is a promising approach to emulate floating-point matrix multiplication on hardware equipped with fast mixed-precision integer matrix-multiplication units.
The algorithm splits the input matrices into integer slices, which are then manipulated using a combination of integer and floating-point arithmetic.
The number of slices is a fundamental design choice: increasing the number of slices increases accuracy at the price of reduced performance.
Our error analysis shows that one can use a different number of slices for the matrices being multiplied, and our experiments suggest that in many scenarios this flexibility can lead to better performance without sacrificing any accuracy.

Standard performance benchmarks, such as the High Performance LINPACK (HPL) benchmark, have not yet embraced the use of these emulation approaches, as they wait for a detailed analysis demonstrating that these techniques can ensure binary64 accuracy.
According to our error analysis, the emulation can become highly inaccurate for badly scaled input matrices, even if a large number of slices is used.
Therefore, this approach cannot be considered a full replacement for binary64 arithmetic---a conclusion reinforced by the fact that, in its current incarnations, the Ozaki scheme does not handle correctly special IEEE 754 values such as negative zeros, infinities, and NaNs.
Further research is needed to determine if and how these special values can be effectively addressed.

\section*{Acknowledgements} This work used resources on the Frank cluster at the Performance Research Laboratory of the University of Oregon. We thank Claude-Pierre Jeannerod (École normale supérieure de Lyon), for pointers to the literature on integer simulation of scalar floating-point arithmetic, and Harun Bayraktar and John Gunnel (NVIDIA), for insightful discussions and early access to the B200~GPUs.

\bibliographystyle{siamplain}
\bibliography{references}

\end{document}